\documentclass[12pt,a4paper]{article}

\usepackage{amsfonts,amsmath,amsxtra,amsthm,amssymb,latexsym}

\newcommand{\ack}{\section*{Acknowledgments}}

\newtheorem{thm}{Theorem}[section]
\newtheorem{cor}[thm]{Corollary}
\newtheorem{lem}[thm]{Lemma}
\newtheorem{prop}[thm]{Proposition}

\newtheorem{defn}{Definition}[section]
\newtheorem{rem}{Remark}[section]

\numberwithin{equation}{section}

\newcommand{\re}{\mathop{\mathrm{Re}}}
\newcommand{\im}{\mathop{\mathrm{Im}}}

\newcommand{\sgn}{\mathrm{sgn}\ }
\newcommand{\dom}{\mathrm{dom}}
\newcommand{\ran}{\mathrm{ran}}
\newcommand{\ess}{\mathrm{ess}}
\newcommand{\disc}{\mathrm{disc}}
\newcommand{\ac}{\mathrm{ac}}

\def\Rset{\mathbb R}
\def\Cset{\mathbb C}
\def\Zset{\mathbb Z}
\def\Nset{\mathbb N}

\def\J{\mathcal{J}}

\def\K{\mathcal{K}}
\def\Ainf{A^{\infty}}
\def\AR{A^0}
\def\ep{\varepsilon}

\newcommand{\Rs}{\mathcal{R}}

\def\A0pm{A_{0\pm}}

\def\mur{\mu^r}
\def\mul{\mu^l}
\def\Cg{c} 
\def\Dg{b} 
\def\A{\mathcal{A}}

\def\D{\mathfrak{D}}
\def\Dmin{\mathfrak{D}_{\min}}

\def\PolP{\mathfrak{p}}
\def\PolQ{\mathfrak{q}}
\def\T{\mathcal{T}}
\def\loc{\mathrm{loc}}

\def\Sec{Sec.}
\def\QQEEDD{}

\title{The similarity problem for $J$-nonnegative
 Sturm-Liouville operators}

\author{Illya M. Karabash,
Aleksey S. Kostenko, and Mark M. Malamud}

\date{}

\begin{document}

\maketitle

\begin{abstract}
Sufficient conditions for the similarity of the
operator $A := \frac {1}{r(x)} \left( -\frac{d^2}{dx^2} +q(x)
\right)$ with an indefinite weight $r(x)=(\sgn x)|r(x)|$
are obtained. These
conditions are formulated in terms of Titchmarsh-Weyl $m$-coefficients. 
Sufficient conditions for the regularity of the critical points
$0$ and $\infty$ of $J$-nonnegative Sturm-Liouville operators are
also obtained. This result is exploited to prove the regularity of
$0$ for various classes of Sturm-Liouville operators. This implies the similarity of the
considered operators to self-adjoint ones. In particular, in the
case $r(x)=\sgn x$ and $q\in L^1(\Rset, (1+|x|)dx)$, we prove that $A$
is similar to a self-adjoint operator if and only if $A$ is
$J$-nonnegative. The latter condition on $q$ is sharp, i.e.,
we construct $q\in \cap_{\gamma <1} L^1(\Rset,
(1+|x|)^\gamma dx)$ such that $A$ is $J$-nonnegative
with the singular critical point $0$. Hence $A$ is not similar to
a self-adjoint operator.
For periodic and infinite-zone potentials,
we show that $J$-positivity is sufficient for the similarity of
$A$ to a self-adjoint operator.  In the case $q\equiv 0$, we prove
the regularity of the critical point $0$ for a wide class
of weights $r$. This yields new results for
"forward-backward" diffusion equations.
\end{abstract}

\quad

similarity,
spectral function of $J$-nonnegative operators, critical points

\quad

 \noindent {\bf Subject classification (MSC2000): } 
 47E05, 34B24, 34B09 (Primary) 34L10, 47B50 (Secondary)\\

\quad

\section{ Introduction}  \label{intro}

Consider the Sturm-Liouville equation
\begin{equation}\label{I_01}
-y''(x)+q(x)y(x)=\lambda\ r(x)y(x),\quad x\in \Rset,
\end{equation}
with a real potential $q \in L^1_{\loc}(\Rset)$ and an
indefinite weight $r\in L^1_{\loc}(\Rset)$.
We assume that $ |r(x)|>0$ a.e. on $\Rset$ and $r$
has only one turning point $x=0$, i.e.,
$ r(x)=(\sgn x)|r(x)|$.

Consider the operator
$L=\frac{1}{|r(x)|}\left(-\frac{d^2}{dx^2}+q(x)\right)$ defined on
its maximal domain $\D$ in the Hilbert space $L^2(\Rset, |r| dx)$. If
$L=L^*$ ($L\geq 0$), then the operator
\begin{equation}\label{I_02}
A := \frac {(\sgn x)}{|r(x)|} \left( -\frac{d^2}{dx^2} +q(x)
\right),\qquad \dom(A)=\D,
\end{equation}
associated with (\ref{I_01}) is called {\em $J$-self-adjoint}
(resp., \emph{$J$-nonnegative}). This means that $A$ is
self-adjoint (nonnegative) with respect to the indefinite inner
product
$[f,g]:=(Jf,g)= \int_\Rset f\overline{g}\, r \, dx , $
where the operator $J$ is defined by
\begin{equation} \label{e J}
(Jf)(x) = (\sgn x) f (x) , \qquad f \in L^2 (\Rset, |r(x)|dx).
\end{equation}

In this paper, we will always assume that $L=L^*$, \ i.e.,
\begin{equation} \label{hyp_limpoint}  \quad \mathit{
the \ differential \ expression \ (\ref{I_01})\ is \ limit \ point
\ at} \ +\infty \ and \ -\infty.
\end{equation}
So the operator $A$ is $J$-self-adjoint. However, it is easy to
see that $A$ is non-self-adjoint in $L^2(\Rset, |r| dx)$ (see
Subsection \ref{sec_II_DifOp}).

The main problem we are concerned with is
\emph{the similarity} of a $J$-nonnegative operator (\ref{I_02})
to a self-adjoint operator. Recall that two closed operators $T_1$
and $T_2$ in a Hilbert space $\mathfrak{H}$ are called \emph{similar} if
there exist a bounded operator $S$ with the bounded inverse
$S^{-1}$ in $\mathfrak{H}$ such that $S\dom(T_1) = \dom(T_2)$ and $T_2 = S
T_1 S^{-1}$.

Ordinary and partial differential operators with indefinite
weights have intensively been investigated during the last two decades
(see \cite{KKLZ84}, \cite{B85}, \cite{CurLan}, \cite{Pyat89},
\cite{Sh93}, \cite{CN95}, \cite{FL96}, \cite{F96}, \cite{Vol96},
\cite{FN98}, \cite{KarMFAT00}, \cite{BV01}, \cite{FSh1}, \cite{FSh2},
\cite{Par03}, \cite{Kar_Mal}, \cite{Zet05}, \cite{Kos05},
\cite{KarKos06}, \cite{BTrunk07}, \cite{KM06} and references
therein).

The similarity of the operator $A$ to a
self-adjoint one is essential for the theory of forward-backward
parabolic equations arising in certain physical models and in
the theory of random processes (see \cite{FK80}, \cite{B85},
\cite{GvdMP87}, \cite{GGvdM88}, \cite{Cur_00}, \cite{FVYu_03},
\cite{Kar_07} and references therein). Theorem
\ref{th_IV.02} of this paper yields new results for
"forward-backward" diffusion equations (see e.g. \cite[Section 5.3]{Kar_07}).

Spectral theory of $\mathcal{J}$-nonnegative operators was
developed by M.G.~Krein and H.~Langer \cite{IKL82,Lan82}
(see Subsection \ref{sec_II_indef}). If the resolvent set $\rho(\mathcal{A})$ of a
$\mathcal{J}$-nonnegative  operator $\mathcal{A}$ is nonempty,
then the spectrum
$\sigma(\mathcal{A})$ of $\mathcal{A}$ is real. Moreover,
$\mathcal{A}$ has a \emph{spectral function} $E_\mathcal{A}
(\cdot)$
with properties similar to that of a spectral function of a
self-adjoint operator. The main difference is the occurrence of
\emph{critical points}. Significantly different behavior of the
spectral function $E_{\mathcal{A}} (\cdot)$ occurs at a \emph{singular
critical point} in any neighborhood of which
$E_{\mathcal{A}} (\cdot)$ is unbounded. A critical point is  \emph{regular} if it is not singular. It should be stressed that only $0$ and
$\infty$ may be critical points of $\mathcal{J}$-nonnegative
operators. Furthermore, $\mathcal{A}$ is similar to a
self-adjoint operator if and only if $0$ and $\infty$ are not singular (see Proposition \ref{p cr=sim}).

If the operator $A$ has a discrete spectrum, the similarity of $A$ to
a self-adjoint operator is equivalent to the \emph{Riesz basis
property} of eigenvectors. For this case, R.~Beals \cite{B85}
showed that the eigenfunctions of Sturm-Liouville problems of
type (\ref{I_01}) form a Riesz basis if $r(x)$ behaves like $(\sgn
x) |x|^\beta$, $\beta >-1/2$, at $x=0$. Improved versions of
Beals' condition were provided in
\cite{CurLan,Pyat89,Sh93,Vol96,F96,Par03}. In \cite{CurLan,F96},
differential operators with nonempty
essential spectrum were considered and the regularity of the
critical point $\infty$ was proved for a wide class of indefinite weight functions. For $J$-nonnegative operators of the form (\ref{I_02}), the result of B.~\'Curgus and H.~Langer \cite[Section
3]{CurLan} is formulated in Proposition \ref{p infReg}.
In particular, it implies the regularity of $\infty$ if
there exist constants $\delta>0$, $\beta_\pm>-1$, and positive functions $\ p_+ \in C^1 [0,\delta]$, $\ p_-\in C^1 [-\delta,0]$
such that
\begin{equation} \label{e BC Simple}
r(x)= (\sgn x) p_\pm (x)|x|^{\beta_\pm} , \quad \pm x \in
(0,\delta).
\end{equation}

The existence of Sturm-Liouville operators of type (\ref{I_02}) with the
singular critical point $\infty$ was established by H. Volkmer
\cite{Vol96} in 1996.  Corresponding examples were constructed
later (see \cite{F96,AbPyat97,F98,Par03,BC04} and references therein).

It turned out that the question of regularity of $0$ is more complicated. Several
abstract similarity criteria
may be found in \cite{Ves72}, \cite{Akop_80}, \cite{Cur_85},
\cite{vCas}, \cite{NabCr}, \cite{MMMCr}, 
\cite{Kap01P}, but it is not easy to apply them to operators
of the form (\ref{I_02}).
First results of this type
were obtained for the operators $(\sgn x) |x|^{-\alpha}
\frac{d^2}{dx^2}$, $\ \alpha>-1$, by B.~\'Curgus, B.~Najman, and
A. Fleige (see \cite{CN95} for the case $\alpha=0$, and
\cite{FN98} for arbitrary $\alpha>-1$). Their approach was
based on the abstract regularity criterion \cite[Theorem
3.2]{Cur_85}. Another approach based on the resolvent criterion of
similarity (see Theorem \ref{t SimCr}) was used by the authors of
the present paper \cite{KarKr98,KarMFAT00,Kar_Mal,Kos05,KM06} as
well as by M.M.~Faddeev and R.G.~Shterenberg \cite{FSh1,FSh2}.
Namely, in \cite{KarKr98,KarMFAT00}, the result of \cite{CN95} was
reproved (see also \cite{Kap01P}). It was shown in \cite{FSh1}
that \emph{if $r(x)=\sgn x$, $\int_{{\Rset}}(1+x^2)|q(x)|dx<\infty$
and $\sigma(A)\subset\Rset$, then $A$ is similar to a self-adjoint
operator}. The case when $q \equiv 0$ and $r(x)\approx
\pm|x|^{\alpha_\pm}, \ \alpha_\pm>-1$, as $\ x\to \pm\infty$, was
considered in \cite{FSh2,Kos05}. A complete analysis for the
case of a finite-zone potential was done in \cite{KM06}.

Our main aim is to present a simple and efficient regularity condition for the critical point $0$ of   operator (\ref{I_02})
and then to apply it to various classes of potentials
(decaying, periodic, and quasi-periodic) as well as to the case when $r(\cdot)$ is nontrivial.
In particular, we show that restrictions
imposed
in \cite{FSh1,FSh2} are superfluous
(see Remarks \ref{rem_first}, \ref{rem_th_FSh}) and give simple
proofs for \cite[Theorem 2.7]{FN98} and
\cite[Corollary 7.4]{KM06}.

Our method is based on two ideas of
\cite{KM06, KarKos06}.
Namely, the resolvent criterion (Theorem \ref{t SimCr})
was used in \cite{KM06} to reduce the
similarity problem to a two weight norm inequality for the
Hilbert transform and to obtain 
similarity conditions in terms of Titchmarsh-Weyl
$m$-coefficients. In particular, \cite[Theorem
5.9]{KM06} states that $A$ is similar to a  self-adjoint operator if
\begin{equation}\label{I_03A}
\sup_{\lambda\in{\Cset}_+}\left|\frac{M_+ (\lambda ) + M_- (\lambda)}{ M_+ (\lambda ) - M_- (\lambda)}\right|<\infty ,\qquad  
\end{equation}
where $M_\pm(\lambda) $ are the Titchmarsh-Weyl $m$-coefficients associated with 
(\ref{I_02}) on $\Rset_\pm$ (explicit definitions are given in
Section \ref{sec_II_indef}).

In this paper we show that
a weaker
form of (\ref{I_03A}) (see Theorem \ref{th_III.1})
remains still sufficient for similarity, and obtain also
its local version using the Krein space approach of \cite{KarKos06}. Namely, if the operator $A$ is $J$-nonnegative and
\begin{equation}\label{I_03}
\sup_{\lambda\in\Omega_R^0}\left|\frac{M_+ (\lambda ) + M_-
(\lambda) - \Cg}{ M_+ (\lambda ) - M_- (\lambda)}\right|<\infty
,\qquad \Omega_R^0:=\{\lambda \in \Cset_+ \, : \, |\lambda| < R \},
\end{equation}
for certain constants $R>0$ and $\ \Cg \in \Rset$, then $0$ is not
a singular critical point of $A$.
%
Combining conditions (\ref{I_03}) and Proposition \ref{p infReg},
we obtain all (sufficient) similarity results of this paper.
However the verification of (\ref{I_03}) requires deep analysis of the $m$-coefficients.

Condition (\ref{I_03A}) is not
necessary \cite[Remark 8.1]{KM06}. Generally, it is violated
for operators considered in Sections \ref{s InfZone} and \ref{Sec_dec}, thought (\ref{I_03}) can be applied
(we do not know whether (\ref{I_03}) is necessary).
Note that the spectral analysis of the finite-zone case \cite[Theorem 7.2]{KM06} was based on the
 similarity criterion (Theorem \ref{t SimCr}) and Muckenhoupt weights rather than on condition
(\ref{I_03A}). The proof of \cite[Theorem 7.2]{KM06} does not require J-nonnegativity of operators, but it is
quite complicated and it is difficult to extend this proof to the operators considered in Sections \ref{s InfZone} and \ref{sec_IV}.

It was proved in \cite{KarKos06} that a condition slightly  weaker than (\ref{I_03}) is necessary for the similarity. Also,  its local version was given (see Theorem \ref{th_III.03}). This result was used to show that the critical
point $0$  of operator $A$ may be singular even 
if $q=0$ (a corresponding example was constructed).
On the other hand, it was proved that there exists a continuous potential
$q\in L^2(\Rset)$ such that the operator $(\sgn x)(-d^2/dx^2+q)$ is
$J$-nonnegative and $0$ is its singular critical point. \emph{The second
aim of this paper is to present an explicit potential with
the above property} (see Theorem \ref{th_VI_1}).

The paper is organized as follows.

In Section \ref{prelim}, we collect necessary definitions and
statements from the spectral theory of Sturm-Liouville operators
and from the spectral theory of $\mathcal{J}$-nonnegative
operators in Krein spaces.

The local regularity condition (\ref{I_03}) is obtained in Section \ref{sec_III}.

In Section \ref{Sec_dec}, we investigate the $J$-self-adjoint operator $A$
with $r(x) =\sgn x$ and $q$ satisfying
\begin{equation} \label{e first}
\int_\Rset (1+|x|)|q(x)|dx<\infty.
\end{equation}
For such operators, we obtain the following criterion.

\begin{thm}\label{t first}
Let $A$ be an operator of the form $(\sgn x)(-d^2/dx^2+q(x))$. If the potential $q$ satisfies (\ref{e first}), then the following statements are equivalent:
\begin{description}
\item $(i)$ $A$ is similar to a self-adjoint operator,
\item $(ii)$ $A$ is $J$-nonnegative (i.e., $L \geq0$),
\item $(iii)$ the spectrum of $A$ is real.
\end{description}
\end{thm}

Under condition (\ref{e first}), $\sigma(L) \cap (-\infty,0)$ may be nonempty but is finite. For this case, we
provide a complete spectral analysis of the operator $A$. Namely,
it is shown that $\sigma_{\ess} (A) = \Rset$, $A$ has no real
eigenvalues, and the discrete spectrum $\sigma_{\disc} (A)$
consists of a finite number of nonreal eigenvalues; we use results of \cite{CurLan} and \cite{KM06} to describe their algebraic and
geometric multiplicities both in terms of definitizing polynomials
and in terms of Titchmarsh-Weyl $m$-coefficients (see Proposition
\ref{th_first_2}).

In Section \ref{Sec_counter}, it is shown that Theorem \ref{t
first} is sharp in the sense that condition (\ref{e first}) cannot
be weaken to $q \in L^1(\Rset, (1+|x|)^{\gamma} dx)$ with $\gamma<1$.
Actually, we construct a potential $q_0$ such that
\begin{description}
\item $(i)$\ \
 $q_0(x)\approx 2(1+|x|)^{-2}$ \emph{as} $|x|\to\infty$,
\item $(ii)$  \ \emph{the operator} $A=(\sgn x)(-d^2/dx^2+q_0(x))$
\emph{is $J$-nonnegative}, \item $(iii)$ \ $0$ \emph{is a singular
critical point of} $A$.
\end{description}
Note that if $r(x) = (\sgn x) |r(x)|$, the regularity of the critical point
$\infty$ of a $J$-nonnegative operator of the form (\ref{I_02}) depends only on local behavior of the weight $r$ in a
neighborhood of $x=0$ (see \cite[Theorem 4.1]{Pyat05}). It appears that the
latter is not true for the critical point $0$. We show  that the
regularity of the critical point $0$ depends not only on
behavior of the weight $r$ at $\infty$ (see \cite[Example 1]{KarKosFAaA07})
but also on local behavior of the potential $q$. This gives an
answer to a one question posed by B.~\'Curgus (see Subsection
\ref{sec_V_II(cur)}).

In Section \ref{s InfZone}, condition (\ref{I_03}) is applied to
operators with periodic potentials.

\begin{thm}\label{th_Per}
Assume that the potential $q \in L^1_{\loc} (\Rset)$ is $\T$-periodic,
$q(x+\T)=q(x)$ a.e., $\T >0$. If the operator $L=-d^2/dx^2 +
q(x)$ is nonnegative, then the operator $A=(\sgn x)L$ is similar
to a self-adjoint operator.
\end{thm}

This theorem can easily be extended to a more general class of
Sturm-Liouville operators with periodic coefficients (see Remark \ref{r gen per}).
Also, a similar result is obtained for the class of infinite-zone
potentials. This class includes smooth periodic potentials.
Generally, infinite-zone potentials are almost-periodic \cite{Lev84}.
For $J$-nonnegative operators with finite-zone potentials, the similarity to a self-adjoint operator was obtained in \cite[Corollary 7.4]{KM06}. We present a simple proof for this result
(see Subsection \ref{ss InfZone}).

In Section \ref{sec_IV}, 
the following theorem is proved.

\begin{thm}\label{th_IV.02}
Let $q\equiv 0$ and $r(x)= \pm p(x)|x|^{\alpha_\pm}$, $x \in \Rset_\pm$,
where $\alpha_\pm > -1$ are constants and the function $p$ is
positive a.e. on $\Rset$. Assume also that
\begin{equation}\label{IV_03}
\pm\int_{\pm1}^{\pm\infty}|x|^{\alpha_\pm/2}|p(x)-c_{\pm}|
dx<\infty,
\end{equation}
with certain constants $c_{\pm}>0$. Then:
\begin{description}
\item $(i)$ $0$ is a regular critical point of the operator $A=-\frac{(\sgn x)}{|r(x)|}\frac{d^2}{dx^2}$;
\item $(ii)$ if the weight $r$ also satisfies the assumptions of Proposition \ref{p infReg}~(i), then the operator $A$ is similar to a self-adjoint one.
\end{description}
\end{thm}

Note that the results of A.~Fleige, B.~Najman \cite[Theorem
2.7]{FN98} and M.M.~Faddeev, R.G.~Shterenberg \cite[Theorem 3]{FSh2}
are particular cases of Theorem \ref{th_IV.02}.
Moreover, we give a short proof of \cite[Theorem
2.7]{FN98}.

Some results of the present paper were announced without proofs in
brief communications \cite{Kar_Mal,Kos06}.
Preliminary version of
this paper was published as a preprint \cite{KAM_07}.

\textbf{Notation:}
Throughout the paper $C_1$, $C_2$, \dots will denote constants that may change from line to line
but will remain independent of the appropriate quantities.
Let $T$ be a linear operator in a Hilbert space $\mathfrak{H}$. In what follows, $\dom (T)$, $\ker (T)$, $\ran (T)$ are the domain,
kernel, range of $T$, respectively; $\sigma(T)$ and $\rho (T)$
denotes the spectrum and the resolvent set of $T$; $ R_T \left(
\lambda \right):=\left( T-\lambda I\right)^{-1} $, $\lambda \in
\rho(T)$, is the resolvent of $T$; $\sigma_p (T)$ stands for the
set of eigenvalues of $T$; the discrete spectrum
$\sigma_{\disc}(T)$ is the set of isolated eigenvalues of finite
algebraic multiplicity;
$\sigma_{\ess}(T):=\sigma(T)\setminus\sigma_{\disc}(T)$ is the
essential spectrum of $T$.

We put $\Cset\pm := \{\lambda \in \Cset :\ \pm \im \lambda
>0\}$, $\Zset_+:=\Nset\cup\{0\}$, $\Rset_+:=[0, +\infty)$, $
\Rset_-:=(-\infty, 0]$. Denote by $\chi_{\mathcal{S}}(\cdot)$ the
indicator function of a set $\mathcal{S}\subset\Rset$,  and
$\chi_\pm(t):=\chi_{\Rset_\pm }(t)$. We write $f\in L^1_{\loc}(\Rset)
(\in AC_{\loc}(\Rset))$ if the function $f$ is Lebesgue integrable
(absolutely continuous) on every bounded interval in $\Rset$;
$f(x)\asymp g(x)\ \ (x\to x_0)$ if both $f/g$ and $g/f$ are
bounded functions in a certain neighborhood of $x_0$; $f(x)\approx
g(x)\ \ (x\to x_0)$ means that $\lim_{x\to x_0} f(x)/g(x)=1$. We
write $f(x)=O(g(x)) \ \left( f(x)=o(g(x)) \right)$ as $x\to x_0$ if
$f(x)=h(x)g(x)$ and $h(x)$ is bounded in a certain neighborhood of
$x_0$ (resp., $\lim_{x\to x_0} h(x)=0$).

 \section{Preliminaries}\label{prelim}

\subsection{Differential operators.}\label{sec_II_DifOp}

Consider the differential expressions
\begin{equation}\label{II_1_01}
\ell [y] := \frac 1{|r|} \left( -y'' +qy \right) \qquad \mathrm{and}
\qquad a [y] :=  \frac 1{r} \left( -y'' +qy \right),
\end{equation}
assuming that $q, r \in L_{\loc}^1(\Rset)$ and $xr(x)>0$ for a.a.  $x\in\Rset$. Let $\D$ be the maximal linear manifold in $L^2
(\Rset, |r(x)|dx)$ on which $\ell
[\cdot]$ and $a[\cdot]$ have a natural meaning:
\begin{equation}\label{II_1_02}
\D:=\{f \in L^2 (\Rset, |r(x)|dx):\ f, f'\in AC_{\loc}(\Rset), \ \ell [f]
\in L^2 (\Rset, |r(x)|dx)\}.
\end{equation}
Define the operators $L$ and $A$ by
\[
\dom(L) = \dom(A) = \D, 
\qquad L f = \ell [f]  \quad \mathrm{and} \quad  A f =  a[f].
\]
The operators $A$ and $L$ are closed in $L^2 (\Rset, |r(x)|dx)$.
In the sequel, (\ref{hyp_limpoint}) is supposed, i. e., $L=L^*$.
It is clear that $A=JL$, where $J^*=J^{-1}=J$ is defined by (\ref{e J}).
Thus, the operator $A$ is $J$-self-adjoint. But $A$ is non-self-adjoint since $A^*=AJ$ and $\dom (A^*) =  J \D \neq  \dom(A)$.

It is obvious that the following restrictions of the operators $L$
and $A$
\begin{eqnarray}
L_{\min}:=L\upharpoonright \Dmin,\quad\qquad A_{\min}:&=&A\upharpoonright
\Dmin,\nonumber\\
 \Dmin:&=&\{ f \in \D\ : f(0)=f'(0) =
0 \},\label{e Amin=}
\end{eqnarray}
are closed densely defined symmetric operators
with equal deficiency indices $(2,2)$.
By $\Dmin^*$ we denote the domain of the adjoint operator
$L_{\min}^*$ of $L_{\min}$. Note that $\Dmin = \D \cap J\D $.
This implies
$\dom(A_{\min}^*)=\dom(L_{\min}^*)=\Dmin^*$ and $\ A_{\min} = JL_{\min}$
(see e.g. \cite{KM06}).
The extensions $\ A_0$ and $\ L_0$
defined by
\begin{eqnarray}
A_0:=A_{\min}^*\upharpoonright \D_0 ,\quad\qquad L_0:&=&L_{\min}^*\upharpoonright \D_0 , \nonumber\\
\D_0 :&=&\{f\in \Dmin^*: f'(+0)=f'(-0)=0\},\label{a 0}
\end{eqnarray}
are self-adjoint operators and $\ A_0=JL_0=L_0 J$.

\subsection{Titchmarsh-Weyl $m$-coefficients.}\label{sec_II_WT}

Let $c(x, \lambda)$ and $s(x, \lambda)$ denote solutions of the initial-value problems
\begin{eqnarray}\label{e eq|r|}
-y''(x)&+&q(x)y(x)=\lambda \ |r(x)|y(x),\qquad x\in \Rset, \\
c(0, \lambda)&=& s'(0, \lambda)=1; \qquad c^{\prime}(0, \lambda)=
s(0, \lambda)=0.\label{II_cs}
\end{eqnarray}
Since equation (\ref{e eq|r|}) is limit-point at $+\infty$,
there exists a unique holomorphic function $m_+(\cdot):\Cset\setminus\Rset\to\Cset$, such that
the solution $s(x, \lambda) - m_+ (\lambda) c(x, \lambda) $
belongs to $L^2(\Rset_+, |r(x)|dx)$ (see e.g. \cite{Tit}). Similarly, the limit point case
at $-\infty$ yields the fact that there exists a unique
holomorphic function $m_-(\cdot):\Cset\setminus\Rset\to\Cset$ such that $
s(x, \lambda) + m_- (\lambda) c(x, \lambda)\in L^2(\Rset_-,
|r(x)|dx)$.
If $\lambda\in\Cset\setminus\Rset$ and $f_\pm(\cdot,\lambda)$
are nontrivial $L^2(\Rset_\pm, |r|dx)$-solutions  of equation (\ref{e eq|r|})
(which are unique up to a multiplicative constant), then
\begin{equation}\label{def_wf}
m_+(\lambda)=-\frac{f_+(+0,\lambda)}{f_+'(+0,\lambda)},\qquad
m_-(\lambda)=\frac{f_-(-0,\lambda)}{f_-'(-0,\lambda)},\qquad
\lambda\notin\Rset.
\end{equation}

The functions $f_\pm(\cdot, \lambda)$ and $m_\pm(\cdot)$  are
called \emph{the Weyl solutions} and \emph{the  Titch-marsh-Weyl
$m$-coefficients (or Titchmarsh-Weyl functions) for} (\ref{e eq|r|})
on $\Rset_\pm$, respectively. We put
\begin{equation}\label{e def psi}
M_\pm (\lambda) := \pm m_\pm (\pm \lambda);\qquad \psi_\pm(x,
\lambda)=(s(x, \pm\lambda)-M_\pm(\lambda)c(x,
\pm\lambda))\chi_\pm(x).
\end{equation}
It is easily seen that $a[\psi_\pm(x,\lambda)]=\lambda
\psi_\pm(x,\lambda)$, where $a[\cdot]$ is defined by
(\ref{II_1_01}). By definition of $m_\pm$,
$\psi_\pm(\cdot, \lambda)\in L^2 (\Rset, |r(x)| dx)$ for all $\lambda
\in \Cset \setminus \Rset$.
The functions
$M_\pm (\cdot)$ are said to be \emph{the
Titchmarsh-Weyl $m$-coefficients for equation (\ref{I_01}) on} $\Rset_\pm$ (associated with the Neumann boundary condition $y'(\pm 0) = 0$).

It is known (see e.g. \cite{Tit}) that  the functions
$\psi_\pm$ and $M_\pm$ are connected by
\begin{equation}\label{217}
\int_{\Rset_\pm} | \psi_\pm(x, \lambda) |^2  |r(x)| dx =  \frac {\im
M_\pm (\lambda)}{\im \lambda} \quad \mathrm{for all} \quad \lambda\in\Cset\setminus\Rset.
\end{equation}
This implies that $M_+$ and $M_-$ (as well
as $m_+$ and $m_-$) belong to the class $(R)$, i.e., \emph{they
are holomorphic in $\Cset\setminus\Rset$, $M_\pm (\overline{\lambda})=
\overline{M_\pm(\lambda)}$, and $\im\lambda \cdot \im M_\pm
(\lambda) \geq 0$  for $ \lambda \in \Cset_+ \cup \Cset_-.$}


\begin{defn}[\cite{KK1}]
An $R$-function $M$ belongs to \\
 $(i)$ the Krein--Stieltjes class $(S)$  if $M$ is holomorphic on $\Cset\setminus\Rset_+$ and $M(\lambda)\geq 0$ for $\lambda<0$;\\
$(ii)$ the Krein--Stieltjes class $(S^{-1})$ if $M$ is holomorphic
on $\Cset\setminus\Rset_+$ and $M(\lambda)\leq 0$ for $\lambda<0$.
\end{defn}
If $M\in(S)$ then it admits the integral representation (see
\cite[\Sec 5]{KK1})
\[
M(\lambda)=c+\int_{-0}^{+\infty}\frac{d\tau(s)}{s-\lambda}\ ,
\qquad \mathrm{where} \quad c\geq 0,\qquad
\int_0^{+\infty}(1+s)^{-1}d\tau(s)<+\infty,
\]
and $\tau:\Rset_+\to\Rset_+$ is a nondecreasing function. This
representation yields that an $S$-function
$M$ is increasing on $(-\infty, 0)$, and
$M(\lambda_0)=0$ for certain $\lambda_0<0$ exactly when $M\equiv
0$. Note also that $M\in (S^{-1})$ if and only if
$\left(-1/M\right)\in (S)$.

The nonnegativity of the self-adjoint operator $L$
can be described in terms of the $m$-coefficients $m_\pm$.

\begin{prop}\label{propII_Sclass}
The operator $L=\frac{1}{|r|}\left(-\frac{d^2}{dx^2}+q\right)$ is nonnegative 
if and only if
$(-1/m_+ - 1/m_-)\in(S^{-1})$.
If, in addition,
$r(x)=-r(-x)$ and $q(x)=q(-x)$ for a.a. $x \in \Rset$,
then $L \geq 0$ exactly when $m_+\in(S)$.
\end{prop}

\begin{proof}
Let $L_{\pm}^D$ be the self-adjoint operators (in $L^2 (\Rset_\pm, |r(x)|dx)$) associated with the Dirichlet boundary value problems
\begin{equation}\label{e LpmD}
-y''(x)+q(x)y(x)=\lambda\ |r(x)|y(x),\qquad x \in \Rset_\pm;\qquad
y(\pm 0)=0.
\end{equation}
Recall that the functions
$\widetilde{m}_\pm :=-1/m_\pm $ are the
Titchmarsh-Weyl $m$-coeffici-ents associated with the problems (\ref{e LpmD}). In particular, $\widetilde{m}_\pm \in (R)$ and
\begin{equation} \label{e spLD}
c(x, \lambda) \pm \widetilde m_\pm (\lambda) s(x, \lambda) \in
L^2(\Rset_\pm, |r(x)|dx) \quad
\mathrm{whenever} \quad \lambda \in \rho(L_{\pm}^D).
\end{equation}
If $L \geq 0 $, then its symmetric restriction $L_{\min}$ defined by
(\ref{e Amin=}) is nonnegative too. Moreover, the extension
$L_+^D \oplus L_-^D$ of $L_{\min}$ corresponding to the
Dirichlet boundary condition at $0$ is a Friedrichs extension, i.e.,
a maximal nonnegative self-adjoint extension of $L_{\min}$
(see \cite{krein_53a} and \cite[Proposition 4]{DM91}).
So $L_{\min} \geq 0$ if and only if $L_{+}^D \geq 0$ and $L_{-}^D \geq 0$. This implies that both $\widetilde{m}_+$ and
$\widetilde{m}_-$ are analytic on $\Cset\setminus\Rset_+$ and
real on $(-\infty, 0)$. It follows from (\ref{e spLD}) that
\[
\{\lambda<0:
\widetilde{m}_+(\lambda)= -\widetilde{m}_-(\lambda)\}=\sigma_p(L)\cap
(-\infty, 0) = \sigma (L)\cap (-\infty, 0).
\]
Since $L\geq0$, we see that $\widetilde{m}_+(\lambda)+\widetilde{m}_-(\lambda)\neq 0$ if
$\lambda<0$. Moreover,
$\widetilde{m}_+(\lambda)+\widetilde{m}_-(\lambda)< 0$ for
$\lambda<0$ since $\widetilde{m}_\pm(-\infty)=-\infty$. 
Thus,
$\widetilde{m}_+ +\widetilde{m}_- \in (S^{-1})$.

If $q$ and $|r|$ are even, then $m_-(\cdot)=m_+(\cdot)$. Hence,
$\widetilde{m}_+ +\widetilde{m}_- =2\widetilde{m}_+(\cdot)=-2/m_+\in
(S^{-1})$ or, equivalently, $m_+\in (S)$. \QQEEDD
\end{proof}

\begin{rem}
In the recent paper \cite{BTrunk07}, the number of negative
squares of self-adjoint operators in Krein spaces were
investigated in terms of abstract Weyl functions
(cf. \cite[Theorem 2]{Kar_Mal}). In particular, Proposition
\ref{propII_Sclass} was proved
 under additional assumptions
 (see Proposition 4.4 and Theorem 4.7 in \cite{BTrunk07}).
      \end{rem}

\subsection{Spectral functions of $J$-nonnegative operators.}\label{sec_II_indef}

Let $\mathfrak{H}$ be a Hilbert space with a scalar product $(\cdot, \cdot )_\mathfrak{H}$. Let $\mathfrak{H}_+$ and $\mathfrak{H}_-$ be closed subspaces of $\mathfrak{H}$ such that $\mathfrak{H} = \mathfrak{H}_+  \oplus \mathfrak{H}_- $. Denote by $P_\pm$ the
orthogonal projections from $\mathfrak{H}$ onto $\mathfrak{H}_\pm$. Put $ \J=P_+ -
P_-$ and $[ \cdot, \cdot ] := (\J \cdot ,\cdot )_\mathfrak{H} $. Then the
pair $\K = (\mathfrak{H} , [\cdot, \cdot])$ is called a \emph{Krein
space} (see e.g. \cite{Lan82,AzJ89} for the original definition). The
form $[ \cdot, \cdot ]$ is called \emph{an inner product} in the
Krein space $\K$ and the operator $\J$ is called \emph{a
fundamental
symmetry}.

Let $T$ be a densely defined operator in $\mathfrak{H}$. By $T^{[*]}$ denote the adjoint of $T$ with respect to $[ \cdot, \cdot ]$.
The
operator $T$ is called \emph{$\J$-self-adjoint}
(\emph{$\J$-nonnegative}) if $ T = T^{[*]} \ $ (resp., $[Tf,f] \geq 0 \ $
for $\ f \in \dom (T)$). It is easy to see that $T^{[*]} := \J T^*
\J$ and $T$ is $\J$-self-adjoint ($\J$-nonnegative)
if and only if $\J T$ is self-adjoint (resp., nonnegative).

Let $\mathfrak{S}$ be the semiring consisting of all bounded
intervals with endpoints different from $0$ and $\pm \infty$ and
their complements in $\overline \Rset := \Rset \cup{ \infty }$.

\begin{thm}[\cite{Lan82}] \label{thII.4.1}
Let $T$ be a $\J$-nonnegative $\J$-self-adjoint operator in
$\mathfrak{H}$ such that $\rho (T) \neq \emptyset$. Then  $\sigma(T) \subset \Rset$ and there exist a mapping $\Delta \rightarrow E(\Delta)$
from $\mathfrak{S}$ into the set of bounded linear operators in
$\mathfrak{H}$ such that the following properties hold ($\Delta ,\Delta' \in
\mathfrak{S}$):
\begin{description}
 \item[(E1)] $E(\Delta \cap \Delta' ) = E(\Delta)
E(\Delta' ) $, \quad $E(\emptyset ) = 0, \quad E(\overline{\Rset}) =
I $, \quad $E(\Delta) = E(\Delta)^{[*]}$;
\item[(E2)] $E(\Delta
\cup \Delta' ) = E(\Delta) + E(\Delta' )$ \quad if \quad $\Delta
\cap \Delta' = \emptyset$;
\item[(E3)] the form $ \pm [
\cdot,\cdot] $ is positive definite on $E (\Delta) \mathfrak{H} $
if $\Delta \subset \Rset_\pm$;
\item[(E4)] $E(\Delta) $ is in the
double commutant of the resolvent $R_T(\lambda)$ and $\sigma (T
\upharpoonright E(\Delta) \mathfrak{H}) \subset
\overline{\Delta}$;
\item[(E5)] if $\Delta$ is bounded, then
$E(\Delta) \mathfrak{H} \subset \dom(T)$ and $T\upharpoonright E
(\Delta) \mathfrak{H} $ is a bounded operator.
\end{description}
\end{thm}

According to \cite[Proposition II.4.2]{Lan82}, $s \in
\{0,\infty\}$ is called \emph{a critical point} of $T$ if, for each $\Delta \in \mathfrak{S}$ such that $s\in \Delta$, the form
$[\cdot,\cdot ]$ is indefinite on $E(\Delta)\mathfrak{H}$ (the latter means that there exist $h_\pm \in E(\Delta)\mathfrak{H}$ such that $\pm [h_\pm,h_\pm ]>0$).
The set of critical points is denoted by $c(T)$. If $\alpha \not
\in c (T)$, then for arbitrary $\lambda_0, \lambda_1 \in \Rset
\setminus c(T)$, $\lambda_0 < \alpha$, $\lambda_1>\alpha$, the
limits
$ \lim_{\lambda \uparrow \alpha} E([\lambda_0,\lambda])$ and
$ \lim_{\lambda \downarrow \alpha} E([\lambda,\lambda_1])$
exist in the strong operator topology; here in the case $\alpha =
\infty$,  $\lambda_1 > \alpha$ ($\lambda \downarrow \alpha$)
means $\lambda_1 > -\infty$ ($\lambda \downarrow -\infty$).
If $\alpha \in c(T)$ and the above limits do still
exist, then $\alpha$ is called \emph{regular critical point} of
$T$, otherwise $\alpha$ is called \emph{singular}.

The following proposition is well known (cf.  \cite[\Sec 6]{Lan82}).
\begin{prop} \label{p cr=sim}
Let T be a $\J$-nonnegative and $\J$-self-adjoint operator in the
Hilbert space $\mathfrak{H}$. Assume that  $\rho(T) \neq
\emptyset$ and $\ker T =\ker T^2$ (i.e., $0$ is either a
semisimple eigenvalue or a regular point of T). Then the following
assertions are equivalent:
\begin{description}
\item $(i)$ T is similar to a self-adjoint operator, \item $(ii)$
$0$ and $\infty$ are not singular critical points of $T$.
\end{description}
\end{prop}

\begin{prop}[\cite{CurLan}, see also \cite{KarKos06}] \label{p RealS}
If the ($J$-self-adjoint) operator $A$ defined by (\ref{I_02}) is
$J$-nonnegative, then its spectrum $\sigma (A)$ is real.
\end{prop}

So any $J$-nonnegative operator of type (\ref{I_02}) has
a spectral function $E_A(\cdot)$.
Note that $\infty$ is always a critical point of $A$,
and $0$ may be its critical point.

\begin{prop}[\cite{CurLan}] \label{p infReg}
Assume that the ($J$-self-adjoint) operator $A$ defined by (\ref{I_02}) is $J$-nonnegative.

\begin{description}
\item $(i)$ Assume that
there exist intervals $ \ \mathcal{I}_\delta^+ = (0,\delta]$,
$ \ \mathcal{I}_\delta^- = [-\delta,0) $, $ \delta >0 $, and constants $s_\pm >0$, $ s_\pm\neq 1$, such that
$r(x) \in AC_{\loc} (\mathcal{I}_\delta^- \cup \mathcal{I}_\delta^-) $, $\left(\frac{r(x)}{r(s_\pm x)}\right)^\prime \in
L^\infty(\mathcal{I}_\delta^\pm)$, and there exist (finite) limits $\lim_{ x\to \pm0}\frac{r(x)}{r(s_\pm x)}\neq s_\pm $.
Then $\infty$ is a regular critical point of $A$.

\item $(ii)$ If $L\geq \varepsilon>0$ and the assumptions of statement (i) are satisfied, then $A$ is similar to a self-adjoint operator.
\end{description}
\end{prop}

Proposition \ref{p infReg} and the slightly stronger condition (\ref{e BC Simple}) follows directly from \cite[Theorem
3.6 (i)]{CurLan}, \cite[Lemma 3.5 (iii)]{CurLan}, and
the remarks in the last two paragraphs of \cite[Subsection 3.2]{CurLan}.

\section{Sufficient conditions for regularity of critical points}\label{sec_III}

Let $A$, $L$, $J$, and $A_{\min}$ be the operators defined in Subsection
\ref{sec_II_DifOp}, and $M_+$, $M_-$ the Titchmarsh-Weyl
$m$-coefficients for (\ref{I_01}) (see Subsection
\ref{sec_II_WT}).

\textbf{3.1} \ \ Our approach to the similarity problem is  based on the resolvent
similarity criterion obtained in \cite{NabCr,MMMCr} (a resolvent similarity criterion, somewhat different from the one given below, was obtained in \cite{vCas}).

     \begin{thm}[\cite{NabCr,MMMCr}] \label{t SimCr}
A closed operator $T$ in a Hilbert space $\mathfrak{H}$ is similar
to a self-adjoint operator if and only if $\sigma (T) \subset \Rset$
and the inequalities
       \begin{eqnarray}\label{ine1}
\sup_{\varepsilon >0}  & \varepsilon   \int_{\Rset} \left\|
\Rs_T
( \eta +i\ep ) f \right\|^2 d\eta \leq K_1 \left\| f\right\| ^2,\quad f \in \mathfrak{H},\\
 \sup_{\varepsilon
>0} & \varepsilon   \int_{\Rset} \left\| \Rs_{T^{*}} ( \eta
+i\varepsilon ) f\right\|^2 d\eta \leq K_{1*} \left\| f\right\|
^2, \quad f \in \mathfrak{H} \label{ine*},
\end{eqnarray}
hold with constants $K_1$ and $K_{1*}$ independent of $f.$
         \end{thm}
\begin{rem} \label{r SimCr}
If $\mathcal{J}=\mathcal{J}^*=\mathcal{J}^{-1}$ and $T$ is a $\mathcal{J}$-self-adjoint operator, then $T^* =
\mathcal{J}T\mathcal{J}$. So (\ref{ine*})  is equivalent to (\ref{ine1}) since in this case
$\|\Rs_{T^{*}}(\lambda)f\|=\|\Rs_{T}(\lambda)f\|$ for all
$\ f\in\mathfrak{H}$, $\ \lambda\in\rho(T)$.
\end{rem}

\textbf{3.2} \ \ For constants $\Dg,\Cg \in \Rset$,
consider the operator $\A_{\Dg,\Cg}:=A_{\min}^*\upharpoonright
\dom(\A_{\Dg,\Cg})$,
\begin{eqnarray} \label{e ACD}
\dom(\A_{\Dg,\Cg})= \left\{f\in \dom(A_{\min}^*) : \
f(+0)-f(-0) =\Cg f'(-0),\right.\nonumber \\
\left.f'(+0)=   \Dg f'(-0), \right\}.
\end{eqnarray}
The operator $A$ defined by
(\ref{I_02}) coincides with $\A_{1,0} $.
Note also that the formal differential expression $ \frac 1r (-\frac{d^2}{dx^2}+q(x)+\Cg\delta'(x))$, where $\delta$
is the Dirac function, may be associated with the operator $A_{1,\Cg}$ (see e.g. \cite{AGHH,Kos05}).

\begin{prop} \label{p Res}
\begin{description}
\item $(i)$ $\A_{\Dg,\Cg}=\A_{\Dg,\Cg}^*$ if and
only if $\Dg=-1$ and $\Cg \in \Rset$.

\item $(ii)$ $\quad \sigma (\A_{\Dg,\Cg}) \setminus \Rset = \{\lambda \in
\Cset_+ \cup \Cset_- : M_- (\lambda) - \Dg M_+ (\lambda) - \Cg = 0\}$.

\item $(iii)$ If $\ \lambda \not \in \Rset \ $ and $\ \lambda \in
\rho(\A_{\Dg,\Cg}) $, then for all $f\in L^2(\Rset, |r|dx)$,
\begin{equation}\label{e III_01}
(\A_{\Dg,\Cg}-\lambda )^{-1}f =(A_0 - \lambda )^{-1} f + \
\frac{\mathcal{F}_- (f,\lambda)-\mathcal{F}_+(f, \lambda)} {M_-
(\lambda) - \Dg M_+ (\lambda) - \Cg} \ \left( \, \Dg \psi_+(\cdot,
\lambda) + \psi_-(\cdot, \lambda) \, \right),
\end{equation}
where
$\mathcal{F}_\pm(f, \lambda):=
\int_{\Rset_\pm}f(x)\psi_\pm(x, \lambda)|r(x)|dx$.
%
\end{description}
\end{prop}

\begin{proof}
\  $(i)$ can be obtained by direct calculation. On the other hand, it follows from the proof of \cite[Proposition 5.8]{KM06}.
Indeed, for the operator $\A_{\Dg,\Cg}$, the matrix $B$ defined by
\cite[formula (5.24)]{KM06} equals $
\left(\begin{array}{cc}
0&\Dg\\
-1&\Cg
\end{array}\right) $, and $\A_{\Dg,\Cg}=\A_{\Dg,\Cg}^*$ exactly when $B=B^*$. 
The proofs of $(ii)$-$(iii)$ are similar to
that of \cite[Proposition
4.4]{KM06} (see also \cite[Lemma 4.1]{KarKos06}). \QQEEDD
\end{proof}

\begin{thm}\label{th_III.1}
If there exists a constant $\Cg \in \Rset$ such that the function
         \begin{equation}\label{e III_03}
\frac{|M_+ (\lambda ) + M_- (\lambda) -
\Cg|}{| M_+ (\lambda ) - M_- (\lambda)|}
        \end{equation}
is bounded on $\Cset_+$, then the operator $A$ is similar to a self-adjoint one.
\end{thm}

\begin{proof}
The proof is similar to that of \cite[Theorem 5.9]{KM06}. We
present a sketch.

Let $\Cg \in \Rset$. Note that $A_0=A_0^*$ (see (\ref{a 0})) and
$\A_{-1,\Cg}=\A_{-1,\Cg}^*$. Hence inequality (\ref{ine1}) holds
for the resolvents of both the operators $A_0$ and $\A_{-1,\Cg}$.
Therefore (\ref{e III_01}) implies that for $f\in L^2(\Rset, |r|dx)$
\begin{equation}\label{e III_3_3}
\sup_{\varepsilon>0} \, \varepsilon
\int_{\mu=-\infty}^{+\infty}
\left\|\frac{\mathcal{F}_{\pm}(f,\mu+i\varepsilon) \, \psi_{\pm}(x, \mu+i\varepsilon)
}{M_+(\mu+i\varepsilon)+M_-(\mu+i\varepsilon)-\Cg}
\right\|_{L^2(|r|dx)}^2
d\mu
\leq C_1\|f\|^2.
\end{equation}
The same arguments and Remark \ref{r SimCr} show that the operator $A=A_{1,0}$ is similar
to a self-adjoint one exactly when
\begin{equation}\label{e III_3_2}
\sup_{\varepsilon>0} \, \varepsilon \int_{\mu=-\infty}^{+\infty}
\left\|\frac{ \mathcal{F}_{\pm}(f,\mu+i\varepsilon) \psi_{\pm}(x, \mu+i\varepsilon) }
{ M_+(\mu+i\varepsilon)-M_-(\mu+i\varepsilon) } \right\|_{L^2(|r|dx)}^2
d\mu
\leq C_2\|f\|^2.
\end{equation}
Combining (\ref{e III_3_3}) with the assumption of the theorem, we get (\ref{e III_3_2}). \QQEEDD
\end{proof}

Theorem \ref{th_III.1} is valid for $J$-self-adjoint (not
necessary $J$-nonnegative) operators of the form (\ref{I_02}). If $\Cg =
0$, this result coincides with \cite[Theorem 5.9]{KM06}.



\begin{thm}\label{th_III.02}
Assume that the operator $A$
is $J$-nonnegative.
If ratio (\ref{e III_03}) is bounded on the set $\Omega_R^0:=\{\lambda \in \Cset_+ \, : \, |\lambda| < R \}$
(on the set $\Omega_R^\infty:=\{\lambda \in \Cset_+ \, : \, |\lambda| > R \}$) for certain constants $R>0$ and $\ \Cg \in \Rset$,
then the point $0$ (resp., the point $\infty$) is not a singular critical point of $A$.
\end{thm}

\begin{proof}
By Proposition \ref{p RealS} and Theorem \ref{thII.4.1}, $A$ has a spectral
function $E_A(\Delta)$. Therefore
$\mathcal{P}_R:=E_A([-R/2, R/2])$ is a bounded $J$-orthogonal
projection. Using properties \textbf{(E1)},\textbf{(E2)}, and \textbf{(E4)} of $E_A(\Delta)$,
we obtain the decomposition
\[
A = \AR \dot + \Ainf , \qquad \AR := A \upharpoonright  \mathfrak{H}_0,
\quad \Ainf := A \upharpoonright  \mathfrak{H}_{\infty}, \quad L^2 (\Rset,
|r|dx) = \mathfrak{H}_0 \dot + \mathfrak{H}_{\infty},
\]
where $ \ \mathfrak{H}_0:= \ran \left( \mathcal{P}_R \right) \ $
and $\ \mathfrak{H}_{\infty}:= \ran \left( I-\mathcal{P}_R \right). $
 Moreover,
\[
\sigma (\AR) \subset [-R/2,R/2], \qquad \sigma(\Ainf) \subset
(-\infty,-R/2] \cup [R/2,+\infty). 
\]
Obviously, $\AR$ is $J$-self-adjoint $J$-nonnegative operator. Note that $\AR$ has the singular critical point $0$ if and
only if so does $A$.

Let us prove that the resolvent of $\AR$ satisfies (\ref{ine1})
if the function
(\ref{e III_03}) is bounded on $\Omega_R^0$.
Indeed, using the last assumption, formula (\ref{e III_01}), and
arguing as in proof of Theorem \ref{th_III.1}, we obtain
\begin{equation} \label{e Ie1}
\varepsilon\int_{\mathcal{I}_\varepsilon}\|(\AR
-(\mu+i\varepsilon))^{-1}f\|^2d\mu=\varepsilon\int_{\mathcal{I}_\varepsilon}\|(A
-(\mu+i\varepsilon))^{-1}f\|^2d\mu\leq C_1 \|f\|^2,
\end{equation}
where
$\mathcal{I}_\varepsilon:=[-\sqrt{R^2-\varepsilon^2},\sqrt{R^2-\varepsilon^2}]$
if $\varepsilon <R$, and $\mathcal{I}_\varepsilon = \emptyset $ if
$\varepsilon \geq R$.
Further, \textbf{(E5)} yields that $\AR$ is bounded. From this and  $\sigma(\AR)\subset [-R/2,
R/2]$, one gets $\|(\AR-\lambda)^{-1}\|\leq C_2|\lambda|^{-1}$ for
$|\lambda|>R$. Hence,
\begin{equation} \label{e Ie2}
\varepsilon\int_{\Rset\setminus \mathcal{I}_\varepsilon}\|(\AR
-(\mu+i\varepsilon))^{-1}f\|^2d\mu\leq C_2\|f\|^2
\int_{\Rset\setminus
\mathcal{I}_\varepsilon}\varepsilon|\mu+i\varepsilon|^{-2}d\mu\leq
C_2\pi\|f\|^2.
\end{equation}
Combining (\ref{e Ie1}) and (\ref{e Ie2}) with Remark \ref{r
SimCr}, we see that $\AR$ is similar to a self-adjoint operator.
Thus $0$ is not a singular critical point of $\AR$.
The proof for the case of the critical point $\infty$ is similar. \QQEEDD
\end{proof}
\textbf{3.3} \ \ In Section \ref{Sec_counter}, we will use the following necessary condition for regularity. 
\begin{thm}[\cite{KarKosFAaA07,KarKos06}]\label{th_III.03}
Assume that the operator $A$ is J-non-negative.
\begin{description}
\item $(i)$  If $0$ is not a singular critical point of $A$ and \
$\ker A =\ker A^2$, then
\begin{equation}\label{302}
\sup_{\lambda\in \Omega_R^0}\left|\frac{\im (M_+ (\lambda
)+M_-(\lambda))} {M_+(\lambda
)-M_-(\lambda)}\right|=C_R<\infty,\quad R>0.
\end{equation}
 \item $(ii)$ If \ $\infty$ is not a singular critical point of $A$, then the function in (\ref{302}) is bounded on $\Omega_R^\infty$ for all
$R>0$.
\end{description}
\end{thm}
\begin{rem}
If $\ \re (M_+ (\lambda) + M_- (\lambda)) -c = O\left( \im (M_+
(\lambda) - M_- (\lambda))\right)$ as $\lambda \to 0$, $\lambda
\in \Cset_+$, the necessary conditions of Theorem \ref{th_III.03}
imply the sufficient conditions of Theorem \ref{th_III.02}. The
results of the following sections show that this is the case for
several classes of coefficients.
\end{rem}

\section{Operators with decaying potentials and regular critical point $0$} \label{Sec_dec}

In this section, we consider the operator
 \begin{equation}\label{first_op}
A=(\sgn x)\left(-\frac{d^2}{dx^2}+q(x)\right),\qquad \dom(A)=\D,
\end{equation}
with the potential $q\in L^1(\Rset)$ having a finite first
moment. That is we consider the case when $r(x)=\sgn x$ and $q$
satisfies
(\ref{e first}).

\subsection{The asymptotic behavior of the Titchmarsh-Weyl $m$-coefficient.}

Since $|r|\equiv 1$, equation (\ref{e eq|r|}) becomes
\begin{equation}\label{dec_01}
-y''(x)+q(x)y(x)=\lambda y(x), \qquad x\in\Rset.
\end{equation}
Note that condition (\ref{e first}) implies that (\ref{dec_01}) is limit point at both $+\infty$ and $-\infty$. Let
$c(\cdot,\lambda)$, $s(\cdot,\lambda)$, and $m_\pm(\cdot)$ be the solutions and the Titchmarsh-Weyl $m$-coefficients of
(\ref{dec_01}) defined as in Subsection \ref{sec_II_WT}.
Denote by $\sqrt{z}$, $\ z\in\Cset \setminus \Rset_+$, the branch of the multifunction $z^{1/2}$ with cut along the positive semi-axis $\Rset_+$ singled out by $\sqrt{-1}=i$.

\begin{lem}\label{lem_decr}
Let 
\begin{equation}\label{first}
\int_{\Rset_+} (1+|x|)|q(x)|dx<\infty.
\end{equation}
Let $s(\cdot,0)$ be the solution of (\ref{dec_01}) 
with $\lambda=0$.
\begin{description}
\item $(i)$ If $s(\cdot,0)$ is unbounded on $\Rset_+$, then
for certain constants $a_+>0$ and $b_+\in\Rset$,
\begin{equation}\label{dec_01_1}
m_+(\lambda)= \frac{a_+}{b_+ \, -i\sqrt{\lambda}}[1+o(1)], \qquad
\lambda\to 0,\quad \lambda\in\Cset\setminus \Rset.
\end{equation}
 \item $(ii)$ If $s(\cdot,0)$ is bounded on $\Rset_+$, then
 for a certain constant $k_+>0$,
\begin{equation}\label{dec_01_1b}
m_+(\lambda)= ik_+ \sqrt{\lambda}\ [1+o(1)], \qquad \lambda\to 0,\quad
\lambda\in\Cset \setminus \Rset.
\end{equation}
\end{description}
%
\end{lem}


\begin{proof}
First note that it suffices to prove (\ref{dec_01_1}) and (\ref{dec_01_1b})
for $\lambda \in \Cset_+$ since $m_+$ is an $R$-function and hence $\overline{m_+(\lambda)}=m_+(\overline{\lambda})$.

 $(i)$ In the case $q\in L^1(\Rset_+)$, the Titchmarsh-Weyl $m$-coefficient admits another
representation (see \cite[Chapter V, \S 3]{Tit}), which is distinct from (\ref{def_wf}). Namely,
\begin{eqnarray}\label{dec_08}
&m_+& (\lambda)=\frac{a(\lambda)}{b(\lambda)},\qquad \qquad\lambda\in\Cset_+, \\
a(\lambda) &=&\frac{i}{2\sqrt{\lambda}}+\frac{i}{2\sqrt{\lambda}}\int_0^{+\infty}q(x)e^{i\sqrt{\lambda}x}s(x,\lambda)dx,\nonumber\\
\qquad
b(\lambda) &=&\frac{1}{2}+\frac{i}{2\sqrt{\lambda}}\int_0^{+\infty}q(x)e^{i\sqrt{\lambda}x}c(x,\lambda)dx, \label{dec_09}
\end{eqnarray}
where the functions $a,b$ are analytic in $\Cset_+$.

In order to estimate $c(x, \lambda)$ and $s(x, \lambda)$, we use
transformation operators preserving initial conditions at the point $x=0$.
Indeed, it follows from \cite[formulas (1.2.9)-(1.2.11)]{Mar77} (see also \cite{Lev84})
that $c(x, \lambda)$ and $s(x, \lambda)$ admit the following representations
\begin{eqnarray}
c(x,\lambda)&=&\cos x\sqrt{\lambda} +\int_{-x}^x K(x,t) \cos t\sqrt{\lambda} dt,\label{dec_02}\\
s(x,\lambda)&=&\frac{\sin x\sqrt{\lambda}}{\sqrt{\lambda}}
+\int_{-x}^x K(x,t) \frac{\sin t\sqrt{\lambda}}{\sqrt{\lambda}}
dt,\label{dec_03}
\end{eqnarray}
where the kernel $K(x, t)$ satisfies the estimates
(see \cite[formulas (1.2.20), (1.2.21)]{Mar77} and also  \cite{Lev84})
\begin{eqnarray}\label{dec_05}
|K(x,t)|\leq
\frac{1}{2} w_0 \left(\frac{x+t}{2}\right) e^{w_1(x)-w_1(\frac{x+t}2)-w_1(\frac{x-t}2)},
\quad 0 \leq |t| < x, \\
 w_0(x):=\int_0^x |q(y)|dy, \qquad w_1(x):=\int_0^x w_0(y)dy.
 \label{dec_05_2}
\end{eqnarray}
Under assumption (\ref{e first}), one can simplify (\ref{dec_05}) as
follows
\begin{equation}\label{dec_06}
|K(x,t)|\leq \frac{1}{2}w_0\left(\frac{x+t}{2}\right)e^{\widetilde{w}_1(\frac{x+t}2)},
\qquad \widetilde{w}_1(x):=\int_0^x\int_y^{+\infty} |q(t)|dtdy
\end{equation}
since inequality $(\ref{e first})$  implies
$\widetilde{w}_1(+\infty)=C_0<\infty$. Hence (\ref{e first}),
(\ref{dec_06}), and (\ref{dec_05_2}) implies
$|K(x,t)|\leq C_1<\infty$ for all $0 \leq |t| < x$. Combining
this fact with (\ref{dec_02}) and (\ref{dec_03}), one obtains
\begin{equation}
|c(x,\lambda)|\leq (1+C_1 x) e^{x \left|\im \sqrt{\lambda} \right|},\qquad
|\sqrt{\lambda}\ s(x,\lambda)|\leq (1+C_1 x) e^{x \left| \im
\sqrt{\lambda} \right|},\label{dec_07'}
\end{equation}
for all $x\in \Rset_+$ and $\lambda\in\Cset_+\cup\Rset.$ We also need the
following inequality  (see \cite[formulas (3.1.28'), (3.1.23)]{Mar77})
\begin{equation}
 |s(x,\lambda)| \leq x e^{x|\im \sqrt{\lambda|}} e^{\widetilde{w}_1(x)}
 \leq C_2 x e^{x|\im \sqrt{\lambda|}}, \quad C_2:=e^{C_0},\label{dec_07}
\end{equation}
which holds for all $\ x\in \Rset_+$, $\lambda\in\Cset_+\cup\Rset$, and is better than (\ref{dec_07'}) as $\lambda\to 0$.

Further, we put
\begin{equation}
\widetilde{a}(\lambda)=1+\int_0^{+\infty}q(t)e^{i\sqrt{\lambda}t}s(t,\lambda)dt,\qquad
\widetilde{b}(\lambda)=\int_0^{+\infty}q(t)e^{i\sqrt{\lambda}t}c(t,\lambda)dt,\label{dec_10}
\end{equation}
\begin{equation}
a_+:=\widetilde{a}(0)=1+\int_0^{+\infty}q(t)s(t,0)dt,\quad
b_+:=\widetilde{b}(0)=\int_0^{+\infty}q(t)c(t,0)dt.\label{dec_11}
\end{equation}
By (\ref{e first}), (\ref{dec_07'}), and (\ref{dec_07}), the
integrals in (\ref{dec_10}) and (\ref{dec_11}) exist and are finite
for all $\lambda\in\Cset_+\cup\Rset$. Note also that $a_+, b_+\in \Rset$
since $q(\cdot), c(\cdot,0)$, and $s(\cdot,0)$ are real functions.

Let us show that $a_+=0$ if and only if $s(x, 0)$ is
bounded on $\Rset_+$. Indeed, integrating the equation
\begin{equation}\label{dec_ker}
-y''(x)+q(x)y(x)=0,\qquad x>0,
\end{equation}
and using $s'(0,0)=1$, we get
\[
s'(x,0)=1+\int_0^{x}q(t)s(t,0)dt,\qquad x\geq0.
\]
By (\ref{dec_11}), $a_+=0$ exactly when
$s'(x,0)=o(1)$ as $x\to+\infty$. On the other hand, equation
(\ref{dec_ker}) with $q(\cdot)$ satisfying (\ref{e first}) has two
linearly independent solutions $y_1(x)$ and $y_2(x)$ such that (see \cite[Theorem X.17.1]{Hartman})
\begin{equation}\label{dec_kernel}
y_1(x)\approx 1,\quad y_1'(x)=o(1);\qquad
y_2(x)\approx x,\quad y_2'(x)\approx 1,
\end{equation}
as $x\to+\infty$.
Hence, $s(x,0)=c_1y_1(x)+c_2y_2(x)$. So we
conclude that $s'(x,0)=o(1)$ as $x\to+\infty$ if and only if
$s(\cdot,0)=c_1y_1(\cdot)\in L^{\infty}(\Rset_+)$.

Note that $c(x,\lambda)$ and $s(x,\lambda)$ are entire functions
of $\lambda$ for every $x\in\Rset_+$. Combining this fact with
(\ref{dec_07'}), (\ref{dec_07}), and first Helly's theorem, we
obtain that functions (\ref{dec_10}) are continuous on
$\Cset_+\cup\Rset$. Due to the assumption $s(\cdot, 0)\notin
L^\infty(\Rset_+)$, we have $a_+\neq0$. Therefore,
\begin{equation}\label{dec_12}
a(\lambda)=i\frac{a_+}{2\sqrt{\lambda}}(1+o(1)),\qquad
b(\lambda)=\frac{1}{2}+i\frac{b_+}{2\sqrt{\lambda}}(1+o(1)),
\end{equation}
as $\lambda\to 0$ and (\ref{dec_01_1}) easily follows from
(\ref{dec_08}) and (\ref{dec_12}).

To complete the proof of $(i)$, it remains to note that $a_+>0$
since $m_+ \in (R)$.

$(ii)$ Let the solution $s(x,0)$ be bounded, i.e., $s(\cdot,0)\in
L^\infty(\Rset_+)$.

Under condition (\ref{first}), for every $\lambda$ in the closed upper half plane $\overline{\Cset_+}$ equation (\ref{dec_01}) has a solution (called \emph{the Jost solution}) that admits the representation
by means of a transformation operator preserving asymptotic behavior at infinity
(see \cite[Lemma 3.1.1]{Mar77}, and also \cite[Chapter I, \S 4]{Lev84})
\begin{equation}
e(x,\lambda):=e^{i\sqrt{\lambda}x}+\int_x^{+\infty}\widetilde{K}(x,t)e^{i\sqrt{\lambda}t}dt
\qquad
x>0,\quad\lambda\in\overline{\Cset_+},\label{dec_17}
\end{equation}
where the kernel $\widetilde{K}(x,t)$ satisfies the following
estimates for $x,t\geq0$
\begin{eqnarray}\label{dec_18}
|\widetilde{K}(x,t)|&\leq&\frac{1}{2}\
\widetilde{\omega}_0\left(\frac{x+t}{2}\right)e^{\widetilde{\omega}(x)},\nonumber\\
&\ &\widetilde{\omega}_0(x):=\int_x^\infty |q(t)|dt,\qquad
\widetilde{\omega}(x)=\int_x^\infty \widetilde{\omega}_0(t)dt.
\end{eqnarray}
Note that, $e(x,\lambda)=e^{i\sqrt{\lambda}x}(1+o(1))$ as $x\to+\infty$. In particular, $e(\cdot,\lambda)$ is the Weyl solution of (\ref{dec_01}) if $\lambda\in\Cset_+$. Moreover,
\[
e(x,0)=1+\int_x^\infty \widetilde{K}(x,t)dt
\]
is a nontrivial 
bounded solution of (\ref{dec_ker}). Hence (cf. (\ref{dec_kernel})),
\[
e(x,0)=c_0's(x,0) \quad \mathrm{with} \quad
(0\neq)\ c_0'=-\widetilde{K}(0,0)+\int_0^\infty K_x'(0,t)dt,
\]
and $ e(0,0)=c_0's(0,0)=0$. Therefore (see \cite[formula (3.2.26)]{Mar77}), $e(0,\lambda)$ has the
form
\begin{equation}\label{dec_19}
e(0,\lambda)=i\sqrt{\lambda}\widehat{K}_1(-\sqrt{\lambda}),\qquad
K_1(x)=\int_x^\infty \widetilde{K}(0,t)dt,
\end{equation}
where 
$\widehat{K}_1(\lambda):=\int_0^\infty K_1(t)e^{-i\lambda t}dt$. Moreover, $\widehat{K}_1$ is continuous at zero since $K_1\in L^1(\Rset_+)$, and $\widehat{c}_0:=\widehat{K}_1(0)\neq 0$ (see the remarks after Eqs. (3.2.25) and (3.2.27) in \cite{Mar77}). 
Noting that $e'(0,0)=c_0's'(0,0)=c_0'\neq 0$, and taking into
account (\ref{def_wf}), we arrive at the desired relation
\begin{equation}
m(\lambda)=-\frac{e(0,\lambda)}{e'(0,\lambda)}=-\frac{i\sqrt{\lambda}\widehat{K}_1(0)}{c_0'}(1+o(1))
,\qquad (\Cset_+\ni)\ \lambda\to 0,
\end{equation}
which proves $(ii)$ with $k_+=-\frac{\widehat{K}_1(0)}{c_0'}$. The
inequality $k_+>0$ follows from the inclusion $m_+ \in (R)$. \QQEEDD
\end{proof}
\begin{rem}
Note that, if $q\in L^1(\Rset_-,(1+|x|)dx)$, then the analogous statements are valid for $m_-$ (with certain constants $a_-,k_->0$, and $b_- \in \Rset$ instead of
$a_+,k_+$, and $b_+$, respectively).
\end{rem}
\begin{prop}\label{prop_p=empty}
Let (\ref{e first}) be fulfilled. Then the operator (\ref{first_op}) 
 has no real eigenvalues, i.e., $\sigma_p(A)\cap\Rset=\emptyset$.
\end{prop}

\begin{proof}
By (\ref{dec_kernel}), $\ker L = \{0\}$. But 
$\ker A=\ker (JL)=\ker L=\{0\}$.

Further, let  $\lambda>0$ and $f(x)\in\ker(A-\lambda)$ (the case
$\lambda<0$ is analogous).
  Then $f\in L^2(\Rset)$ solves (\ref{dec_01}) with $\lambda>0$. Under assumption (\ref{e first}), equation (\ref{dec_01}) has two linearly independent solutions of the form (see \cite[Lemma 3.1.1]{Mar77}) 
\begin{eqnarray*}
e_+(x,\lambda)&=&e^{i\sqrt{\lambda}x}+\int_x^{+\infty}\widetilde{K}(x,t)e^{i\sqrt{\lambda}t}dt,\\
e_-(x,\lambda)=e^{-i\sqrt{\lambda}x}&+&\int_x^{+\infty}\widetilde{K}(x,t)e^{-i\sqrt{\lambda}t}dt;\qquad
x\geq0,
\end{eqnarray*}
with $\widetilde{K}$ satisfying (\ref{dec_18}). 
So $f(x)=c_+ e_+ (x,\lambda)+ c_- e_- (x,\lambda)$ for $x>0$ with
certain $c_\pm\in\Cset$. Hence (\ref{dec_18}) implies that $f (x)\approx c_+ e^{
i\sqrt{\lambda}x}+c_- e^{- i\sqrt{\lambda}x}$ as $x\to+\infty$ (see \cite[formula (3.1.20)]{Mar77}).
The latter yields $c_+=c_-=0$ since $f\in L^2(\Rset)$. Therefore,
$f(x)=0$, $x>0$. Since $f$ is a solution of (\ref{dec_01}), we get
$f \equiv 0$.
\QQEEDD
\end{proof}

\begin{rem} \label{rem_FN ext}
Assume that $q$ satisfies (\ref{e first}) on $\Rset_+$ and that the
minimal symmetric operator $L^+_{\min}$ associated with the spectral problem
\[
-y''(x)+q(x)y(x)=\lambda \ y(x),\qquad
x\geq0,\quad y(0)=y'(0)=0,
\]
is nonnegative in $L^2(\Rset_+)$. The Friedrichs (hard) extension $L^D_+=(L^D_+)^*$ of $L_{\min}^+$ is
determined by the Dirichlet boundary condition at zero (for definitions and basic facts on M.G. Krein's extension theory of nonnegative operators see \cite[\Sec 109]{AGII}).
The corresponding $m$-coefficient is
$\widetilde{m}_+(\cdot) \ (=-1/m_+(\cdot))$. Lemma
\ref{lem_decr} shows that $s(\cdot, 0)\in
L^\infty(\Rset_+)$ exactly when $ \widetilde{m}_+(-0)=+\infty$.
It follows from \cite{krein_53a} (see also \cite[Proposition 4]{DM91}) that $ \widetilde{m}_+(-0)=+\infty$ holds if and only if $L^D_+$ is the Krein--von Neumann (soft) extension of $L_{\min}^+$. The latter means that the operator $L_{\min}^+$ has a unique nonnegative self-adjoint extension.
Thus, Lemma \ref{lem_decr} leads to the following criterion
: \emph{$L^D_+$ is
a unique nonnegative self-adjoint extension of the nonnegative
operator $L_{\min}^+$ if and only if $s(\cdot, 0)\in
L^\infty(\Rset_+)$}.
\end{rem}

\subsection{The case of the nonnegative operator $L$.}

The proof of Theorem \ref{t first} is contained in this and the
next subsections. The most substantial part, the implication $(ii)
\Rightarrow (i)$, is given by the following result:

\begin{thm}\label{p first}
Let $A=(\sgn x)(-d^2/dx^2 + q(x))$ and let $q(\cdot)$ satisfy
(\ref{e first}). If the operator $A$ is $J$-nonnegative, then it is
similar to a self-adjoint operator.
\end{thm}

\begin{proof} 
Assume that the operator $A$ is $J$-nonnegative. By Proposition
\ref{p RealS}, $\sigma(A)\subset\Rset$. Proposition \ref{p infReg} implies that
$\infty$ is a regular critical point of $A$.  Moreover,
(\ref{e first}) implies $\ker A=\{0\}$ (see Proposition
\ref{prop_p=empty}). Hence the similarity of $A$ is equivalent to
the nonsingularity of the critical point zero of the operator $A$
(see Proposition \ref{p cr=sim}).

By Lemma \ref{lem_decr} and (\ref{e def psi}), one of the
asymptotic formulas (\ref{dec_01_1}), (\ref{dec_01_1b})  holds for
the function $m_+(\lambda)=M_+(\lambda)$. And the same is true for
$m_- (\lambda) = -M_-(-\lambda)$.
Consider the following four cases.\\
$(\mathrm{a})$ \emph{Let the solution $s(\cdot, 0)$ of (\ref{dec_01}) be
bounded on $\Rset$}, $s(\cdot, 0)\in L^\infty(\Rset)$. By Lemma
\ref{lem_decr} (iii), for $(\Cset_+\ni) \lambda \to 0$ we get
\[
M_+(\lambda)= i k_+\sqrt{\lambda}\ (1+o(1)),\quad M_-(\lambda)=
k_-\sqrt{\lambda}\ (1+o(1)); \quad k_\pm>0.
\]
Therefore, we obtain as $\lambda\to 0$
\[
\frac{M_+(\lambda)+M_-(\lambda)}{M_+(\lambda)-M_-(\lambda)}=\frac{ik_{+}\sqrt{\lambda}+k_{-}\sqrt{\lambda}}{ik_{+}\sqrt{\lambda}-k_{-}\sqrt{\lambda}}(1+o(1))=
\frac{ik_+ + k_-}{ik_+ -k_-}(1+o(1)).
\]
$(\mathrm{b})$ \emph{Let $s(\cdot, 0)\notin L^\infty(\Rset_+)$, but
$s(\cdot,0)\in L^\infty(\Rset_-)$}. Then, by Lemma \ref{lem_decr},
\[
M_+(\lambda)=\frac{a_+}{b_+ -i\sqrt{\lambda}}(1+o(1)),\qquad
M_-(\lambda)= k_-\sqrt{\lambda}(1+o(1)); \qquad \lambda\to 0,
\]
where $a_+>0, \ b_+\in\Rset, \ \mathrm{and}\ k_->0$. Hence
we get
\begin{equation}\label{IV_17'}
\frac{M_+(\lambda)+M_-(\lambda)}{M_+(\lambda)-M_-(\lambda)}=\frac{a_++k_{-}\sqrt{\lambda}(b_+
-i\sqrt{\lambda})}{a_+-k_{-}\sqrt{\lambda}(b_+
-i\sqrt{\lambda})}(1+o(1))=1+O(\sqrt{|\lambda|}).
\end{equation}
$(\mathrm{c})$ The case when \emph{$s(\cdot, 0)\in L^\infty(\Rset_+)$ and
$s(\cdot,0)\notin L^\infty(\Rset_-)$} is similar to (b).

$(\mathrm{d})$ \emph{Let  $s(\cdot, 0)\notin L^\infty(\Rset_+)$ and $s(\cdot,
0)\notin L^\infty(\Rset_-)$}. Then, by Lemma \ref{lem_decr} (ii), one
gets as $\lambda\to 0$
\[
M_+(\lambda)=a_+(b_+-i\sqrt{\lambda})^{-1}(1+o(1)),\quad
M_-(\lambda)=-a_-(b_- +\sqrt{\lambda})^{-1}(1+o(1)),
\]
where $ a_\pm>0$ and $b_\pm\in\Rset$. Hence,
\begin{eqnarray}
\frac{M_+(\lambda)+M_-(\lambda)-\Cg}{M_+(\lambda)-M_-(\lambda)}  \qquad\qquad\qquad\qquad\qquad\qquad\qquad\qquad\qquad\qquad\nonumber\\
\approx \frac{a_+(b_-+\sqrt{\lambda})-a_-(b_+-i\sqrt{\lambda})-\Cg
(b_+-i\sqrt{\lambda})(b_-+\sqrt{\lambda})}{a_+(b_-+\sqrt{\lambda})+a_-(b_+-i\sqrt{\lambda})} \quad 
\label{dec_20}
\end{eqnarray}
as $\lambda\to 0$. If $b_+\cdot b_-=0$, then the left part of
(\ref{dec_20}) with $\Cg=0$ has the asymptotic behavior similar
to one of the cases (a),(b), or (c). Otherwise, we put
$\Cg:=a_+b_+^{-1}-a_-b_-^{-1}$ and get
\[
\frac{M_+(\lambda)+M_-(\lambda)-\Cg}{M_+(\lambda)-M_-(\lambda)}=\frac{(a_+-\Cg
b_+)\sqrt{\lambda}-(a_-+\Cg
b_-)\sqrt{-\lambda}}{a_+(b_-+\sqrt{\lambda})+a_-(b_++\sqrt{-\lambda})}\left(1+o(1)\right)=O(1)
\]
as $\lambda\to 0.$ From the above considerations, we conclude that
there exists $\Cg\in\Rset$ such that ratio (\ref{e III_03}) is
bounded in a neighborhood of zero. By Theorem \ref{th_III.02},
zero is not a singular critical point of $A$.  Combining this fact
with Propositions \ref{p cr=sim} and \ref{p RealS}, we complete
the proof of the similarity of $A$ to a self-adjoint operator. \QQEEDD
\end{proof}

In passing, we have proved the following fact for any (not
necessarily $J$-nonnegative) operator $(\sgn x)(-\frac{d^2}{dx^2}+q(x))$ with $q$ satisfying (\ref{e first}).
\begin{prop}\label{prop_0reg}
Let $A=(\sgn x)(-\frac{d^2}{dx^2}+q(x))$ with $q$ satisfying
(\ref{e first}). Then there exist $\Cg\in \Rset$ such that ratio
(\ref{e III_03}) is bounded in a neighborhood of zero.
\end{prop}

\begin{rem}\label{rem_definitizable}
 It should be pointed out that if $L$ is nonnegative, then only the case (d) in the proof of Theorem \ref{p first} can be realized. Actually, if $s(\cdot,0)\in L^\infty(\Rset_+)$, then  (\ref{dec_01_1b}) yields $(-m_+(x))^{-1} \uparrow +\infty$ as $x\uparrow -0$. Therefore $(-m_+)^{-1}+(-m_-)^{-1}$ takes positive values on $\Rset_-$. But $L\geq 0$ and Proposition \ref{propII_Sclass}
implies $(-m_+)^{-1}+(-m_-)^{-1}\in (S^{-1})$. This contradiction
shows that $s(\cdot, 0)\notin L^\infty(\Rset_\pm)$.
\end{rem}

\subsection{The operator $L$ with negative eigenvalues.}

It is known that
under condition (\ref{e first}), the negative
spectrum of the operator $L=JA=-d^2/dx^2 + q(x)$ consists of at
most finite number $\kappa_-(L)$ of simple eigenvalues  and (see \cite[Theorem 5.3]{Ber_Shub})
\[
\kappa_-(L)\leq 1+\int_\Rset |x|q_-(x)dx,\qquad
q_-(x):=(|q(x)|-q(x))/2.
\]
 So Propositions 1.1 and 2.5 of \cite{CurLan} imply that $A$ is a
\emph{definitizable} operator (for the definitions and basic facts
see \cite{JonL_79, Lan82, CurLan} and \cite[Appendix B]{F96}). The
latter means that $\rho(A) \neq \emptyset$ and there exists a real
polynomial $\PolP$ such that $[\PolP (A) f,f ] \geq 0$ for all $f
\in \dom (A^k)$, where $k=\deg \PolP$; the polynomial $\PolP$ is
called \emph{definitizing}. Since $\kappa_-(JA)$ is finite, there
is a definitizing polynomial $\PolP$ of minimal degree  and of the
form (see \cite[Eq. (1.2)]{CurLan})
\begin{equation} \label{e DefP}
\PolP (z) = z \PolQ (z) \overline{\PolQ(\overline{z})}, \qquad
\deg \PolQ\leq \kappa_-(L).
\end{equation}
The polynomial $\PolQ(z)$ is uniquely determined under the
assumption that it is monic polynomial and all its zeros belongs
to $\Cset_+ \cup \Rset$. A definitizable operator admits a spectral
function $E(\Delta)$ with, possibly, some critical points (which
belong to the set $\infty \cup \{ \lambda \in \Rset : \PolP
(\lambda)=0\}$). The properties of $E (\Delta)$ similar to that of
$E (\Delta)$ from Theorem \ref{thII.4.1}.

B. \'Curgus and H. Langer \cite{CurLan} investigated nonreal
spectrum of indefinite $\J$-self-adjoint ordinary differential
operators $\mathcal{A}$ assuming that $\J \mathcal{A}$ has a finite number of negative
eigenvalues. The following result follows from \cite[Subsection
1.3]{CurLan}.

\begin{prop} \label{p zeroP}
Let $A=(\sgn x)(-d^2/dx^2 + q(x))$ and $q \in L^1 (\Rset,
(1+|x|)dx)$. Let $\PolQ$ be defined by (\ref{e DefP}). Then:
\begin{description}
\item $(i)$ $\lambda \in \Rset\setminus\{0\}$ is a zero of
$\PolQ(\cdot)$ if and only if it is a critical point of $A$; \\in
this case, $\lambda$ is also an eigenvalue of $A$.
\item $(ii)$
$\lambda \in \Cset_+ \ (\Cset_-)$ is a zero of $\PolQ(\cdot)$
(resp., $\PolQ(\overline{\cdot})$) if and only if it is a nonreal
eigenvalue of $A$; in this case, the algebraic multiplicity of
$\lambda$ is finite.
\end{description}
\end{prop}

Taking Propositions \ref{prop_p=empty} and \ref{prop_0reg} into
account, we obtain the following description for essential and
discrete parts of the operator $A$.

\begin{thm}\label{th_first_2}
Let $A=(\sgn x)(-d^2/dx^2 + q(x))$ and $q(\cdot)$ satisfy
(\ref{e first}). Then:
 \begin{description}
 \item $(i)$ The nonreal spectrum $\sigma(A) \setminus \Rset $
 is finite and consists of eigenvalues of finite algebraic multiplicity. If $\lambda_0 \in \Cset \setminus \Rset$ is an eigenvalue of $A$, then its algebraic multiplicity
 is equal to the multiplicity of $\lambda_0$ as a zero of the holomorphic function $M_+ (\lambda) - M_- (\lambda)$. Its geometric multiplicity equals $1$.
 \item $(ii)$ $\sigma_p(A)=\sigma_{\disc} (A)=\sigma(A) \setminus \Rset$,\ \  $\sigma_{\ess} (A) = \Rset$, and there exist a skew direct decomposition
 $L^2(\Rset)=\mathfrak{H}_{\ess}\dotplus\mathfrak{H}_{\disc}$ such that
    \begin{eqnarray}
     A&=&A_{\ess}\dotplus A_{\disc},\nonumber\\
       A_{\ess}=A\upharpoonright(\dom(A)&\cap&\mathfrak{H}_{\ess}),\qquad
      A_{\disc}=A\upharpoonright(\dom(A)\cap\mathfrak{H}_{\disc}),\nonumber\\
     \sigma_{\disc} (A)=\sigma(A_{\disc}) \
     (&=&\sigma(A) \setminus \Rset),
     \qquad \sigma_{\ess}(A)=\sigma(A_{\ess}) \ (=\Rset);\nonumber
    \end{eqnarray}
 the subspace $\mathfrak{H}_{\disc}$ is finite-dimensional.
 \item $(iii)$ $A_{\ess}$ 
 is similar to a self-adjoint operator.
\end{description}
\end{thm}

\begin{proof}
\ $(i)$ follows from Proposition \ref{p zeroP} (i) and
\cite[Proposition 4.3 (5)]{KM06}.

$(ii)$ follows from (i) and Proposition \ref{prop_p=empty} (see
e.g. \cite[Section 6]{KM06}). Note only that
$\sigma_{\ess}(A)=\sigma_{\ess}(A_0)=\Rset$ (see e.g.
\cite[Proposition 4.3 (1)]{KM06}).

$(iii)$ The operator $A$ is a definitizable and admits a spectral
function $E_A (\Delta)$. By Proposition \ref{prop_p=empty}, $\sigma_p (A)\cap \Rset = \emptyset$. So Proposition \ref{p zeroP} (i) implies that $\PolQ$
has no real zeros and that the only possible critical points of
$A$ are zero and infinity (actually, $0$ and $\infty$ {\em are} critical points).  Further, $\infty$ is a regular
critical point due to \cite[Theorem 3.6]{CurLan}. Using
Proposition \ref{prop_0reg} and arguing as in the proof of
Theorem \ref{th_III.02}, one can prove that zero is not a
singular critical point of $A$. Hence $A_{\ess}$, the part of $A$
corresponding to the real spectrum, is similar to a 
self-adjoint operator $T$.
\QQEEDD
\end{proof}

\begin{cor} \label{c CritReal}
Let $A=(\sgn x)(-d^2/dx^2 + q(x))$ and $q(\cdot)$ satisfy
(\ref{e first}). Then:
\begin{description}
\item $(i)$ $\lambda_0$ is an eigenvalue of $A$ if and only if it
is a zero of $\PolQ(z)\overline{\PolQ(\overline{z})}$; moreover,
its algebraic multiplicity coincides with the multiplicity as a
zero of $\PolQ(z)\overline{\PolQ(\overline{z})}$. \item $(ii)$
$\sigma(A)\subset\Rset$ if and only if $A$ is $J$-nonnegative.
\end{description}
\end{cor}

\begin{proof}
\ $(i)$ Since $\sigma_p (A) \cap \Rset = \emptyset$, Proposition \ref{p zeroP} (i) implies $\PolQ(0) \neq
0$. It follows from these
facts that equality holds in \cite[formula (1.3)]{CurLan}.
Combining this and \cite[Proposition II.2.1]{Lan82}, we see that
the degree $\deg \PolP$ of polynomial
$\PolP(z)=z\PolQ(z)\overline{\PolQ(\overline{z})}$ is greater or
equal than $2\kappa_-(JA)$. From this and (\ref{e DefP}), we
obtain $\deg \PolP = 2\kappa_-(JA)+1$ and $\deg \PolQ =
\kappa_-(JA)$. Applying the equality in \cite[formula
(1.3)]{CurLan} and \cite[Proposition II.2.1]{Lan82} again, one
gets statement (i).

$(ii)$ For the case $L \geq 0$, see Proposition \ref{p RealS}. If $A$ is not $J$-nonnegative, then $\kappa_-(JA) \geq 1$
and therefore $\PolQ(\cdot) \not \equiv 1$. So $\PolQ(\cdot)$ has at least one
zero $\lambda_1$, which is an eigenvalue of $A$ due to statement
(i) and is nonreal due to Proposition \ref{prop_p=empty}.
 \QQEEDD
\end{proof}
Now we are ready to prove Theorem \ref{t first}.
\begin{proof}[Proof of Theorem \ref{t first}]
Note that the implication $(ii)\Rightarrow (i)$ follows from
Theorem \ref{p first}. The implication $(i)\Rightarrow (iii)$ is
obvious. To complete the proof it suffices to mention
that the equivalence $(ii)\Leftrightarrow(iii)$ was established in
Corollary \ref{c CritReal} (ii). \QQEEDD
\end{proof}
Recall that the function $M_+ (\cdot) - M_- (\cdot)$ is
holomorphic in $\Cset\setminus\Rset$. The next result follows easily
from Theorem \ref{th_first_2} (i) and Corollary \ref{c CritReal}
(i).
\begin{cor}\label{col_first}
Let $A=(\sgn x)(-d^2/dx^2 + q(x))$ and $q(\cdot)$ satisfy
(\ref{e first}). Assume also that $A$ is not $J$-nonnegative.
\begin{description}
\item $(i)$ Let $ \{ z_j \}_1^n $ be the set of nonreal zeros of
the function $M_+ (\cdot) - M_- (\cdot)$, and let $ \{ k_j \}_1^n
$ be their multiplicities. Then $\PolP = z \prod_1^n
(z-z_j)^{k_j}$ is a definitizing polynomial of minimal degree for
$A$. \item $(ii)$ $A$ is similar to a normal operator if and only
if $k_j=1$ for all $1\leq j \leq n$.
\end{description}
\end{cor}



\begin{rem}\label{rem_first}
Under the additional assumption $q \in L^1(\Rset, (1+|x|^2)dx$, the
equivalence $(i)\Leftrightarrow(iii)$ in Theorem \ref{t first} was
proved in \cite{FSh1} by using another approach.
Note also that inclusion $\sigma(A)\subset\Rset$ was established in \cite[Corollary 4]{FSh1} under the assumption $m_\pm\in (S)$ (cf. 
Proposition \ref{propII_Sclass} of the present paper).
\end{rem}


\section{Operators with decaying potentials and singular critical point $0$}\label{Sec_counter}

If $r(x) = \sgn x$, then the operator $A$ defined by (\ref{I_02}) is similar to a self-adjoint one whenever $L(=JA)$ is uniformly positive (see Proposition \ref{p infReg}). If $0 \in \sigma_{\ess} (L)$, then it may occur that $0$ is a critical point of $A$. Sturm-Liouville operators of type $-\frac{d^2}{r(x)dx^2}$ with the \emph{singular} critical point $0$ were constructed in \cite{KarKos06}. 
A $J$-nonnegative operator of type $(\sgn x)(-d^2/dx^2+q(x))$ with
the singular critical point $0$ have not been constructed, but
existence of such an operator was proved in \cite[Section
6.2]{KarKos06}. The goal of this section is to construct
an explicit example of such type. Our example also shows that
condition (\ref{e first}) in Theorem \ref{t first} cannot be
weaken to $q \in L^1(\Rset, (1+|x|)^{\gamma} dx)$ with $\gamma<1$.


\subsection{Example.}

\begin{lem}\label{lem_VI_1}
Let
\begin{equation}\label{VI_01}
q_0(x)=-\chi_{[0,
\pi/4]}(x)+2\frac{\chi_{(\pi/4,+\infty)}(x)}{(1+x-\pi/4)^2},\qquad
x\in\Rset_+.
\end{equation}
Then the function
\begin{equation}\label{VI_02}
m_0(\lambda)=\frac{\sin(\pi\sqrt{\lambda+1}/4)/\sqrt{\lambda+1}+m_1(\lambda)\cos(\pi\sqrt{\lambda+1}/4)}{\cos(\pi\sqrt{\lambda+1}/4)-m_1(\lambda)\sqrt{\lambda+1}\sin(\pi\sqrt{\lambda+1}/4)},\quad\lambda\in\Cset_+,
\end{equation}
where
\begin{equation}\label{VI_03}
m_1(\lambda)=\frac{1-i\sqrt{\lambda}}{1-i\sqrt{\lambda}-\lambda},\qquad
\lambda\in \Cset_+,
\end{equation}
is the Titchmarsh-Weyl $m$-coefficient of the boundary value
problem
\begin{equation}\label{VI_04}
-y''(x)+q_0(x)y(x)=\lambda y(x),\quad x\geq0;\quad y'(0)=0.
\end{equation}
\end{lem}
\begin{proof} Consider the Sturm-Liouville equation
\begin{equation}\label{VI_06}
-y''(x)+\frac{2}{(1+x)^2}y(x)=\lambda y(x),\quad x\geq0.
\end{equation}
It is easy to check that
$f_1(x,\lambda)=e^{i\sqrt{\lambda}(x+1)}(\sqrt{\lambda}+i/(x+1))$
solves (\ref{VI_06}) and  $f_1(\cdot,\lambda)\in L^2(\Rset)$ for
$\lambda\in\Cset_+$. Further, 
$f_1(0,\lambda)=e^{i\sqrt{\lambda}}(\sqrt{\lambda}+i)$ and
$f_1'(0,\lambda)=e^{i\sqrt{\lambda}}(-\sqrt{\lambda}+i\lambda-i)$.
By (\ref{def_wf}), we get that (\ref{VI_03}) is the
Titchmarsh-Weyl $m$-coefficient of (\ref{VI_06}) associated with  the Neumann boundary condition at zero.

Using (\ref{VI_01}), we obtain that the function
\begin{eqnarray}\label{VI_09}
f_0(x,\lambda)=&\left(f_1(0,\lambda)\cos((x-\frac{\pi}4)\sqrt{\lambda+1})+\right.
\nonumber\\
&\left. f_1'(0,\lambda)\frac{\sin((x-\frac{\pi}4)\sqrt{\lambda+1})}{\sqrt{\lambda+1}}\right)\chi_{[0,\frac{\pi}4]}(x)+\\
&f_1(x-\frac{\pi}4,\lambda)\chi_{(\frac{\pi}4,+\infty)}(x),\qquad x\geq0,\nonumber
\end{eqnarray}
is the Weyl solution of (\ref{VI_04}) for $\lambda\in\Cset_+$. To
complete the proof, it remains to substitute (\ref{VI_09}) in
(\ref{def_wf}). \QQEEDD
\end{proof}
Let us consider the indefinite Sturm-Liouville operator
\begin{equation}\label{VI_05}
A=(\sgn x)\left(-\frac{d^2}{dx^2}+q_0(|x|)\right),\qquad
\dom(A)=W_2^2(\Rset),
\end{equation}
with $q_0$ defined by (\ref{VI_01}).
\begin{thm}\label{th_VI_1}
Let $A$ be the operator defined by (\ref{VI_05}) and
(\ref{VI_01}). Then:
\begin{description}
\item $(i)$ $A$  is $J$-self-adjoint, $J$-nonnegative, and
$\sigma(A)\subset \Rset$.
\item $(ii)$  $0$ is a simple eigenvalue of
$A$, i.e., its algebraic multiplicity is 1.
\item $(iii)$ $0$ is a singular critical point of $A$.
\item $(iv)$ $A$ is not similar to a self-adjoint operator.
\end{description}
\end{thm}

\begin{proof}
\ $(i)$ Note that $q_0$ is bounded on $\Rset$. Hence $A$ is
$J$-self-adjoint. Next, we show that the operator
$L=JA=-d^2/dx^2+q_0(|x|)$  is nonnegative. The potential is even,
hence, by Lemma \ref{lem_VI_1},
$m_+(\lambda)=m_-(\lambda)=m_0(\lambda)$ (see
(\ref{VI_02})). It is easy to see that $m_1$ is a Krein-Stieltjes
function, $m_1\in(S)$, since it is analytic and positive on
$(-\infty, 0)$. It is not difficult to see that the latter
implies $m_0\in (S)$. Proposition \ref{propII_Sclass} yields
$L\geq0$. Hence $A=JL$ is $J$-nonnegative and, by Proposition
\ref{p RealS}, $\sigma(A)\subset \Rset$.

$(ii)$ It is easily seen that $\lim_{\lambda\to 0}\lambda
m_0(\lambda)=k\neq 0$.
So $\lambda=0$ is the
eigenvalue of the problem (\ref{VI_04}). Hence $c(x, 0)\chi_+(x)\in
L^2(\Rset_+)$. Furthermore, $q_0(|x|)$ is even, hence $c(x,
0)\chi_-(x)\in L^2(\Rset_-)$ and $c(x, 0)\in\ker L$. Since
$s(x,0)\notin L^2(\Rset)$, we get $\ker L=\{a\ c(\cdot, 0):
a\in\Cset\}$. The equality $\ker A=\ker L$ implies $0\in
\sigma_p(A)$.

Further, by (\ref{VI_02}) and (\ref{VI_03}), we get
\[
m_0(\lambda)=\frac{1+m_1(\lambda)}{1-m_1(\lambda)}(1+O(\lambda))=\left(1+\frac{2}{\sqrt{-\lambda}}-\frac{2}{\lambda}\right)\left(1+O(|\lambda|)\right),\quad
|\lambda|\to 0.
\]
Note that $M_+(\cdot)=-M_-(-\cdot)=m_0(\cdot)$ since
$m_+(\cdot)=m_-(\cdot)=m_0(\cdot)$. Hence,
\begin{equation}\label{VI_07}
\frac{\im(M_+(iy)+M_-(iy))}{M_+(iy)-M_-(iy)}=\frac{\im\
m_0(iy)}{\re\ m_0(iy)} =
\frac{1/\sqrt{2y}+1/y}{1+1/\sqrt{2y}}\left(1+O(y)\right)\approx\sqrt{\frac{2}{y}}
\end{equation}
as $ y\to+0$. Combining (\ref{e III_01}) with (\ref{217}),
(\ref{VI_07}), and the inequality $\|(A_0-iy)^{-1}\|\leq y^{-1}$,
after simple calculations we arrive at
\[
\|(A-iy)^{-1}\|\leq
O(y^{-3/2}),\qquad y\to +0.
\]
Therefore, $\ker A=\ker A^2$. This completes the proof of $(ii)$.

$(iii)$ Combining (\ref{VI_07}) with Theorem \ref{th_III.03} (i),
we conclude that $0$ is a singular critical point of $A$.

$(iv)$ follows from Proposition \ref{p cr=sim} and $(iii)$.
\QQEEDD
\end{proof}

\subsection{On a question of B. \'Curgus.}\label{sec_V_II(cur)}

It is known that infinity is a critical point of the operator
(\ref{I_02}). Moreover, the results of \cite{CurLan, Vol96, F96,
Par03} shows that the regularity of the critical point $\infty$ of
a definitizable operator of type (\ref{I_02}) depends only on
behavior of the weight function $r$ in a neighborhood of its
turning point (in our case, in a neighborhood of $x=0$). At
$6^{th}$ Workshop on Operator Theory in Krein Spaces (TU Berlin,
2006), B.\'Curgus posed the following problem: \emph{does the
regularity of the critical point zero of a $J$-nonnegative
operator of type (\ref{I_02}) depend only on behavior of the
coefficients $q$ and $r$ at infinity}?

Below we give the negative answer to this question.

Consider the operator
\[
A_1=(\sgn
x)\left(-\frac{d^2}{dx^2}+2\frac{\chi_{(\pi/4,+\infty)}(|x|)}{(1+|x|-\pi/4)^2}\right),\qquad
\dom(A)=W_2^2(\Rset).
\]
It is easy to see that $A_1$ is $J$-self-adjoint and
$J$-nonnegative since the potential is bounded and positive on
$\Rset$. Arguing as in the proof of Lemma \ref{lem_VI_1}, we obtain
that the corresponding Titchmarsh-Weyl $m$-coefficients are
 \[
M_+(\lambda)=-M_-(-\lambda)=m_2(\lambda):=\frac{\sin(\pi\sqrt{\lambda}/4)/\sqrt{\lambda}+m_1(\lambda)\cos(\pi\sqrt{\lambda}/4)}{\cos(\pi\sqrt{\lambda}/4)-m_1(\lambda)\sqrt{\lambda}\sin(\pi\sqrt{\lambda}/4)},
 \]
where $m_1(\cdot)$ is given by (\ref{VI_03}). Since $m_1(\lambda)=1+O(\sqrt{\lambda})$ as $\lambda\to 0$, we easily get $m_2(\lambda)=(1+\pi/4 +O(\sqrt{\lambda}))$ as $\lambda\to
0$.
Hence we obtain
\[
\lim_{\overline{\Cset_+}\ni\lambda\to
0}\left|\frac{M_+(\lambda)+M_-(\lambda)}{M_+(\lambda)-M_-(\lambda)}\right|=
\left|\frac{(1+\pi/4)-(1+\pi/4)}{(1+\pi/4)+(1+\pi/4)}\right|=0<\infty,
\]
and, by Theorem \ref{th_III.02}, $0$ is not a singular
critical point of $A_1$.

On the other hand, the operator $A$ considered in the previous
subsection
is an additive perturbation of $A_1$ by a potential with a compact
support. However, $0$ is a singular critical point of $A$ due to
Theorem \ref{th_VI_1} (iii). Thus, \emph{the regularity of the
critical point zero of operator $(\ref{I_02})$ depends not
only on behavior of the weight function $r$, but also on local
behavior of the potential $q$.}



\section{Operators with periodic and almost-periodic potentials} \label{s InfZone}

Throughout this section we assume $r (x) =\sgn x$, so the
operators $L$ and $A$ have the forms $L=-d^2/dx^2 +q(x)$ and
$A=(\sgn x) L$. All the asymptotic formulas in this section are
considered in $\Cset_+$.

\subsection{The case of a periodic potential $q$.}
\label{ss Periodic}

First, we consider the case of $\T$-periodic potential $q \in
L^1_{\loc} (\Rset)$, i.e., $q(x+\T) =q(x)$ a.e. on $\Rset$, $\T>0$. It
is known that in this case equation (\ref{e eq|r|}) is limit point
at both $+ \infty$ and $-\infty$. Hence, the maximal operator $L$
corresponding to the differential expression $-d^2/dx^2 + q(x)$ is
self-adjoint in $L^2 (\Rset)$.

Let $c(x,\lambda)$ and $s(x,\lambda)$ be the functions defined by
(\ref{e eq|r|}),  (\ref{II_cs}). Recall that for any $x\in\Rset$,
$c(x,\lambda)$, $s(x,\lambda)$, $c'(x,\lambda)$, and
$s'(x,\lambda)$ are entire functions of $\lambda$, hence so are
\begin{equation} \label{e F}
\Delta_+(\lambda) := \frac{c(\T,\lambda)+s'(\T,\lambda)}2\qquad
\mathrm{and}\qquad \Delta_-(\lambda) :=\frac{
c(\T,\lambda)-s'(\T,\lambda)}2\ .
\end{equation}
The function $2\Delta_+(\cdot)$ is the trace of the monodromy
matrix and it is called Hill's discriminant (or the Lyapunov
function).

As before, we denote by $\widetilde m_\pm (\lambda)$ ($ m_\pm
(\lambda)$) the Titchmarsh-Weyl $m$-coefficient for (\ref{e eq|r|})
on $\Rset_\pm$ corresponding to the Dirichlet (Neumann, resp.)
boundary condition at $0$. Then,
\begin{equation}\label{e ML Per}
\widetilde m_{\pm}(\lambda) = - \frac{1}{m_\pm (\lambda)} =
  \frac {\mp \Delta_-(\lambda) + \sqrt{\Delta_+^2(\lambda)-1}}
{ s(\T,\lambda) } \ ,
       \end{equation}
where the branch  of the multifunction
$\sqrt{\Delta_+^2(\lambda)-1}$ is  chosen  such  that both
$\widetilde m_{\pm}(\cdot)$ (and so $ m_{\pm}(\cdot)$) belong to
the class $(R)$.
For continuous $q(\cdot)$, formula (\ref{e ML
Per}) may be found, e.g., in \cite[\Sec 21.2]{Tit}, the proof of (\ref{e ML Per}) for $q \in L^1 [0,\T]$ is the same.

\begin{lem}\label{l la0}
Let $L$ be a Sturm-Liouville operator with a $\T$-periodic
potential $q\in L^1_{\loc}(\Rset)$. Let also $\lambda_0:=\inf
\sigma(L)$. Then:
\begin{description}
\item $(i)$ \ $(-\infty<)\ \lambda_0$ is a first order zero of
$\Delta_+(\lambda)-1$ and $\Delta_+'(\lambda_0)<0$; \item $(ii)$ \
$s(\T, \lambda_0)>0$.
\end{description}
\end{lem}

This statement is well known for the case
of continuous $q$ (see e.g. \cite[\Sec 21.4]{Tit}). For the case $q \in L^1 [0,\T]$, it can be obtained, e.g., from \cite[Sections 12 and 13]{Weid87}.

\begin{proof}
(ii) follows from \cite[Theorem 13.7 (a)]{Weid87}.

(i) The proofs of Theorems 12.5 (c), 12.7, and 13.10 in \cite{Weid87} show that $\lambda_0$ is the first eigenvalue of the corresponding periodic problem, $\Delta_+(\lambda_0)=1$, $\Delta_+(\lambda) >1$ for $\lambda<\lambda_0$, and
$\Delta_+(\lambda) <1$ for $\lambda-\lambda_0>0$ small enough.
So the order $n_{\lambda_0}$ of $\lambda_0$ as a zero of the entire function $\Delta_+(\lambda)-1$ is an odd number.
Let us show that $n_{\lambda_0} = 1$ (and therefore, $\Delta_+^\prime(\lambda_0)<0$).
It follows from (\ref{e ML Per}) and statement (ii) that $\widetilde m_+(\lambda) \approx C_1 + C_2 (\lambda - \lambda_0)^{n_{\lambda_0}/2}$ as $\lambda \to \lambda_0$, where $C_1$, $C_2$ are real constants and $C_2 \neq 0$.
So if $n_{\lambda_0} \geq 3$, then $\widetilde{m}_+ \not \in (R)$, a contradiction.
 \QQEEDD
\end{proof}



\begin{proof}[Proof of Theorem \ref{th_Per}]
 Consider the operator $L= -d^2/dx^2 +q(x)$ with a $\T$-periodic potential $q$ and assume that $\lambda_0 (= \inf \sigma (L)) \geq 0$.
It follows from (\ref{e def psi}) and (\ref{e ML Per}) that
Titchmarsh-Weyl $m$-coefficients for the operator $A=(\sgn x)L$ have
the form
       \begin{equation}\label{e M Per}
M_\pm (\lambda ) =   \frac { s(\T,\pm\lambda)}
{\Delta_-(\pm\lambda) \mp \sqrt{\Delta_+^2(\pm\lambda)-1} } \ .
      \end{equation}

By Proposition \ref{p infReg}, $\infty$ is a regular critical
point of $A$. At the same time, by Proposition \ref{p infReg}, it
suffices to consider only the case $\lambda_0 = 0 $.

Assuming $\lambda_0 =0$, consider two cases.

$(a)$ Let $\Delta_-(0)=0$.
Lemma \ref{l la0} (i) yields that $\lambda_0=0$ is a first order
zero of the entire function $\Delta_+(\lambda)-1$. By  Lemma
\ref{l la0} (ii), $s(\T,0)>0$ and, therefore, (\ref{e M Per})
implies
\begin{equation}\label{e asM Case1}
M_\pm (\lambda)  =   \frac { s(\T,0)(1+O(\lambda))}
{\pm\lambda(\Delta_-'(0)+O(\lambda)) \mp \sqrt{\pm\lambda
(2\Delta_+'(0)+O(\lambda))} }=\pm\  i\frac { C_1}{\sqrt{ \pm
\lambda} }[1+O(\sqrt{\lambda})],
\end{equation}
as $ \lambda \to 0$, where $C_1=s(\T,0)/\sqrt{-2\Delta_+'(0)}>0$. Substituting
(\ref{e asM Case1}) for $M_\pm(\cdot)$ in (\ref{e III_03}) with
$\Cg=0$, we see that Theorem \ref{th_III.02} implies that
$0$ is not a singular critical point of $A$.

$(b)$ Suppose $\Delta_-(0)\neq 0$. Note that $\Delta_-(\lambda)$
and $\Delta_+(\lambda)$ are real if $\lambda\in\Rset$. Combining
(\ref{e M Per}) with Lemma \ref{l la0} (ii), we get
\begin{equation}\label{e asM Case2}
M_\pm (\lambda) =   \frac { s(\T,0)} {\Delta_-(0) \mp\ i
C_2\sqrt{\pm\lambda } }[1+O(\sqrt{\lambda})],
\qquad \lambda \to 0,
\end{equation}
with $C_2 = \sqrt{-2\Delta_+'(0)}>0$. Using Theorem
\ref{th_III.02} with
$\Cg=2s(\T,0)/\Delta_-(0)\in\Rset\setminus\{0\}$, we see that $0$ is
not a singular critical point of $A$.

Thus the operator $A$ is $J$-nonnegative and has no singular
critical points. Moreover, $\ker A=\ker L$, and $\ker L=\{0\}$ since $q$ is $\T$-periodic
(see e.g. \cite[\Sec 12]{Weid87}).
Proposition \ref{p cr=sim} completes the proof of
Theorem \ref{th_Per}. \QQEEDD
\end{proof}



\begin{rem} \label{r gen per}
Let $\T$-periodic functions $p$, $q$, and $\omega$
be such that $\frac 1p, q, \omega \in L^1_\loc (\Rset)$
and $p, \omega >0$ a.e. on $\Rset$.
Then the operator
\[
(Ly)(x):=
\frac{1}{\omega(x)}\left(-(p(x)y'(x))'+q(x)y(x)\right)
\]
defined on the maximal domain in $L^2 (\Rset, \omega (x)dx )$
is self-adjoint and semi-bounded from below. Moreover, equation (\ref{e ML Per}) and Lemma \ref{l la0} hold with the same proofs.
Therefore the proof of Theorem \ref{th_Per} shows that { \em
$0$ is not a singular critical point of the operator $A:= (\sgn x) L$ whenever $A$ is J-nonnegative.} If additionally the critical point $\infty$ is regular, then {\em $A$ is similar to a self-adjoint operator.}
For instance, the latter holds if $p, \frac 1p \in L^{\infty} (-\delta,\delta)$ for certain $\delta>0$ and the function $r(x):=(\sgn x) \omega (x)$ satisfies the conditions of Proposition
\ref{p infReg} (i) (see \cite[\Sec 3]{CurLan}).
\end{rem}

\subsection{Infinite-zone and finite-zone potentials.}
\label{ss InfZone}

In this subsection we consider the cases of (real) infinite- and
finite-zone potentials.

Following \cite{Lev84}, we briefly recall definitions. First
note that the spectrum of the operator $L=-d^2/dx^2 + q(x)$ with
an infinite-zone potential $q$ is absolutely continuous and has
the zone structure, i.e.,
\begin{equation} \label{e sinf}
\sigma (L) = \sigma_{\ac} (L)=[\mur_0 , \mul_1 ] \cup [\mur_1 ,
\mul_2 ] \cup \cdots,
\end{equation}
where $\{ \mur_j \}_0^\infty$ and $\{ \mul_j \}_{j=1}^{\infty}$
are sequences of real numbers such that
\begin{equation} \label{e mu<mu}
\mur_0 < \mul_1 < \mur_1 < \dots < \mur_{j-1} < \mul_j < \mur_j <
\dots \quad ,
\end{equation}
and
\[
\lim_{j\to \infty} \mur_j = \lim_{j\to \infty} \mul_j = +\infty.
\]
In the case of a finite-zone potential, the corresponding
sequences $\{ \mur_j \}_0^N$, $\{ \mul_j \}_{j=1}^N$ are finite,
$N<\infty$, the spectrum of $L$ is also absolutely continuous and
is given by
\begin{equation} \label{e sfin}
\sigma (L) = \sigma_{\ac} (L)=[\mur_0 , \mul_1 ] \cup [\mur_1 ,
\mul_2 ] \cup \cdots \cup [\mur_N , + \infty) \ .
\end{equation}

Let $N \in \Zset_+ $. Consider also sets of real
numbers
$\{ \xi_j \}_1^N $ and $\{\epsilon_j \}_1^N $ such that 
 $ \xi_j \in [\mul_j , \mur_j ] $ and $\epsilon_j \in \{-1,+1\}$
 for all $j \leq N$.
Define polynomials $R(\lambda)$, $P(\lambda)$, and $Q(\lambda)$ by
     \begin{eqnarray}\label{PolPR}
& P(\lambda ) = \prod_{j=1}^N (\lambda - \xi_j ) , \qquad
R(\lambda ) = (\lambda - \mur_0)
\prod_{j=1}^N (\lambda - \mul_j ) (\lambda- \mur_j) , \\
& Q(\lambda) = P(\lambda)\sum_{j=1}^N \frac{\epsilon_j
\sqrt{-R(\xi_j)}}{P'(\xi_j)(\lambda-\xi_j)}. \label{PolQ}
      \end{eqnarray}
Then  there exists (see \cite[Lemma 8.1.1]{Lev84}) a  real
polynomial $S(\lambda)$ of degree $\deg S =N+1$ such that
\begin{equation}   \label{PolinomS}
S(\lambda) = \prod_{j=0}^N (\lambda - \tau_j ) , \quad \tau_0 \in
( - \infty , \mur_0 ] , \quad \tau_j \in [\mul_j , \mur_j ] ,
\quad j\in \{1, \dots, N\},
\end{equation}
and the following identity  holds
       \begin{equation} \label{e PSQR}
P(\lambda) S(\lambda) - Q^2 (\lambda ) = R(\lambda ) .
     \end{equation}

According to \cite[formulas (8.1.9) and (8.1.10)]{Lev84}
 the functions
     \begin{equation}\label{e ML FZ+}
m_{\pm}(\lambda)  :=  \pm  \frac {P(\lambda)}{ Q(\lambda) \mp i
\sqrt{R( \lambda)}}
       \end{equation}
are the Titchmarsh-Weyl $m$-coefficients 
corresponding to the Neumann boundary value problems on $\Rset_{\pm}$
for some Sturm-Liouville operator  $L=-d^2/dx^2 +q(x)$ with a
\emph{quasi-periodic} potential $q=\bar q$ (see e.g. \cite[\Sec 10.3]{Lev84}). Here  the
multifunction $\sqrt{R( \cdot)}$ is considered on $\Cset$ with cuts
along the union of intervals (\ref{e sfin}).  The branch $\sqrt{R(
\cdot)}$ of the multifunction is  chosen in such a way that
$\sqrt{R(\lambda_0+i0)}>0$ for some $\lambda_0\in (\mur_{N},
+\infty).$ So
both  $m_{\pm}(\cdot)$ belong to the class $(R)$. In this case the
spectrum of $L$ is given by (\ref{e sfin}).

       \begin{defn}[\cite{Lev84}] \label{d FZp}
A real quasi-periodic potential $q$ is called finite-zone if the
Titchmarsh-Weyl $m$-coefficients $m_{\pm}$  admit the
representations  (\ref{e ML FZ+}).
       \end{defn}

Note that if the potential $q$ is $\T$-periodic and the equation
$\Delta_+^2 (\lambda)= 1$ (see (\ref{e F})) has a finite number
of simple roots, then $q$ is a finite-zone potential (see
\cite[Sections 7.4 and 8.1]{Lev84}). Moreover, in this case $\mur_j$ and $\mul_j$
denote simple roots of $\Delta_+^2(\lambda)-1=0$ listed in the
natural order. Note also that every finite-zone potential $q$ is
bounded and its n-th derivative $\frac {d^n}{dx^n}q$ is bounded on
$\Rset$ for any $n\in\Nset$ (see \cite[\Sec 8.3]{Lev84}).


A criterion of the similarity to a self-adjoint
operator for (not necessary $J$-nonnegative) operator $A=(\sgn
x)(-d^2/dx^2+q(x))$ with a finite-zone potential was
obtained in \cite[Theorems 7.1 and 7.2]{KM06}. For the case of a
$J$-nonnegative operator $A$, we present a new simple proof of
\cite[Corollary 7.4]{KM06} based on Theorem \ref{th_III.02}.

\begin{thm}[\cite{KM06}]\label{th_FinZ}
Let $q(x)$ be a finite-zone potential and $\mur_0 \geq 0$. Then
the operator $A=(\sgn x)(-d^2/dx^2 + q(x))$ is similar to a self-adjoint
operator.
\end{thm}
\begin{proof}
Consider the operator $L= -d^2/dx^2 +q(x)$ with a finite-zone
potential $q$ and assume that $L \geq 0$. This is equivalent to
$\mur_0 \geq 0$ due to (\ref{e sfin}).

Combining (\ref{e def psi}) with (\ref{e ML FZ+}) and (\ref{e
PSQR}), we get
       \begin{equation}\label{e M FZ}
M_\pm (\lambda ) = \frac {P(\pm \lambda)}{ Q(\pm \lambda) \mp i
\sqrt{R(\pm \lambda)}} =  \frac {Q(\pm \lambda) \pm i \sqrt{R(\pm
\lambda)}  }{ S (\pm \lambda)} .
      \end{equation}
It is easy to see that
\begin{equation} \label{e MinInf}
M_\pm (\lambda) = \pm \frac {i}{ \sqrt{\pm
\lambda}}[1+O(\lambda^{-1/2})],  \qquad \lambda \to \infty,\qquad
\lambda \in \Cset_+.
\end{equation}
This implies that the function $(M_+ +M_-)(M_+ -M_-)^{-1}$ is
bounded in a certain neighborhood of $\infty$. So $\infty$ is a
regular critical point due to Theorem \ref{th_III.02}.

Let us prove that $0$ is not a singular critical point. As in the
periodic case, we note that $0 $ is not a critical point if
$\mur_0 > 0 $.  Further, assume that $\mur_0 =0$ and consider the
cases analogous to that of the proof of Theorem \ref{th_Per}.

$(a)$ Let $\tau_0 = 0 \ (= \mur_0)$, where $\tau_0$ is 
defined in (\ref{PolinomS}). Then $R(0)=S(0)=0$, and it follows
from (\ref{e PSQR}) that $Q(0)=0$. By definition, $P (0) = P
(\mur_0) \neq 0$ and, therefore, (\ref{e ML FZ+}) implies that
(\ref{e asM Case1}) holds
with $C_1 = \frac{\prod_{j=1}^{N}  \xi_j}{\left(\prod_{j=1}^{N} \mul_j \mur_j \right)^{-1/2}} >0 $.

$(b)$ Let $\tau_0 \neq 0 $ (actually, this yields $\tau_0 <0$, see
(\ref{PolinomS})). Then $S(0) \neq 0$. Further, $R(0)=0$, $P(0)
\neq 0$ and (\ref{e PSQR}) implies that $Q(0) \neq 0$. Using the
second representation of $M_\pm (\lambda)$ from (\ref{e M FZ}),
one can check that
\begin{equation}\label{e asM Case2fz}
M_\pm (\lambda) = C_2 \pm i \, C_3 \sqrt{ \pm \lambda} +
o(|\lambda| ^{1/2}) , \qquad \lambda \to 0 ,
\end{equation}
where $C_2 = Q(0)/S(0) \in \Rset \setminus \{ 0 \}$
 and $C_3 = |C_1/S(0)| >0$.


The arguments of Subsection \ref{ss Periodic} conclude the proof.
\QQEEDD
\end{proof}

In the proof of Theorem \ref{th_FinZ}, we have shown that $\infty$ is a regular critical point of $A$ using the asymptotic formula (\ref{e MinInf}) for $M_\pm$ and the regularity condition, Theorem
\ref{th_III.02}.  On the other hand, this fact follows from
Proposition \ref{p infReg}.



Now consider infinite sequences $\{ \mur_j \}_0^\infty$, $\{
\mul_j \}_{1}^{\infty}$, $\{ \xi_j \}_1^\infty $, and
$\{\epsilon_j \}_1^\infty $ such that $ \xi_j \in [\mul_j , \mur_j
] $, $\epsilon_j \in \{-1,+1\}$
 for all $j \geq 1$, and assumptions (\ref{e mu<mu}) and inequalities
 \begin{eqnarray} \label{e SumGap}
\sum_{j=1}^\infty \mur_j (\mur_j - \mul_j) < \infty, \qquad
\sum_{j=1}^\infty \frac 1{\mul_j} < \infty . \quad
\end{eqnarray}
 are fulfilled.
For every $N \in \Nset$, put
\begin{eqnarray}
 & g_N  =  \prod_{j=1}^N \frac{\xi_j - \lambda }{\mul_j}, \quad
f_N  = (\lambda - \mur_0) \prod_{j=1}^N \frac{\lambda - \mul_j}
{\mul_j}
\frac{\lambda- \mur_j }{\mul_j}, \label{e fNgN} \\
 & k_N (\lambda)  =  g_N (\lambda)
\sum_{j=1}^N \frac{\epsilon_j \sqrt{-f_N
(\xi_j)}}{g'_N(\xi_j)(\lambda-\xi_j)} , \quad h_N (\lambda)  =
\frac{f_N (\lambda) + k_N^2 (\lambda)}{g_N (\lambda)}. \label{e
kNhN}
\end{eqnarray}

It is easy to see from (\ref{e SumGap}) that $g_N$ and $f_N$
converge uniformly on every compact subset of $\Cset$. Denote $
\lim_{N \to \infty} g_N (\lambda) =: g (\lambda)$, $ \lim_{N \to
\infty} f_N (\lambda) =: f (\lambda)$. \cite[Theorem 9.1.1]{Lev84}
states that there exist limits $\lim_{N \to \infty} h_N (\lambda)
=: h (\lambda)$, $ \lim_{N \to \infty} k_N (\lambda) =: k
(\lambda)$ for all $\lambda \in \Cset$. Moreover, the functions $g$,
$f$, $h$, and $k$ are holomorphic in $\Cset$.

It follows from \cite[Subsection 9.1.2]{Lev84} that
 the functions
     \begin{equation}\label{e ML inf}
m_{\pm}(\lambda)  :=  \pm  \frac {g(\lambda)}{ k(\lambda) \mp i
\sqrt{f( \lambda)}}
       \end{equation}
are the Titchmarsh-Weyl $m$-coefficients on  $\Rset_{\pm}$
(corresponding to the Neumann boundary conditions) for some
Sturm-Liouville operator  $L=-d^2/dx^2 +q(x)$ with a real bounded
potential $q(\cdot)$. The branch $\sqrt{f( \cdot)}$ of the
multifunction is  chosen  such  that both  $m_{\pm}(\cdot)$ belong
to the class $(R)$.
          \begin{defn}[\cite{Lev84}] \label{d infZp}
A real potential $q$ is called an infinite-zone potential if the
Titchmarsh-Weyl $m$-coefficients $m_{\pm}$ admit representations
(\ref{e ML inf}).
       \end{defn}
Let $q$ be an infinite-zone potential defined as above. Since $q$
is bounded, the operator $L=-d^2/dx^2 +q(x)$ is self-adjoint. Its
spectrum is given by (\ref{e sinf}). B. Levitan proved that under
the additional condition $\inf
 (\mul_{j+1} - \mul_j) >0$, the potential $q$ is \emph{almost-periodical} (see \cite[Chapter 11]{Lev84}).
 Note that for a $\T$-periodic
 potential $q$ the first inequality in (\ref{e SumGap}) implies $q\in W^2_2[0,\T]$,
 and the second inequality in (\ref{e SumGap}) obviously follows from
 asymptotic formulas for the periodic (anti-periodic) eigenvalues
 (see \cite[\Sec 1.5]{Mar77} for details).





The following theorem is the main result of this subsection.
\begin{thm}\label{th_InfZ}
Let $L=-d^2/dx^2 + q(x)$ be a Sturm-Liouville operator with an
infinite-zone potential $q$. Assume also that the spectrum
$\sigma(L)$ satisfies (\ref{e SumGap}) and $L \geq 0$ (i.e., $\mur_0
\geq 0$). Then the operator $A=(\sgn x)L$ is similar to a
self-adjoint operator.
\end{thm}

The asymptotic formula (\ref{e MinInf}) does not hold true in the
infinite-zone case. Therefore, we use Proposition \ref{p infReg}
to prove that $\infty$ is a regular critical point. The rest of
the proof is also close to subsection \ref{ss Periodic}.

\begin{proof}
It is sufficient to consider the case $\mur_0 = 0$. Recall that
the functions $g$, $f$, $k$, and $h$ defined above are holomorphic
in $\Cset$. Moreover, $g$ and $f$ admit the
following representations
\[ 
g (\lambda)=  \prod_{j=1}^\infty \frac{\xi_j - \lambda }{\mul_j},
\quad f (\lambda)= \lambda \prod_{j=1}^N \frac{\lambda - \mul_j}
{\mul_j} \frac{\lambda- \mur_j }{\mul_j}, 
\]  
where the infinite products converge uniformly on all compact
subsets of $\Cset$ due to assumptions (\ref{e SumGap}) (see
\cite[Section 9]{Lev84}). From this and $\xi_j>\mur_0 =0$,
$j\in \Nset$, we see that
\begin{equation} \label{e f=0 gneq0}
f(0) = 0 , \qquad g(0) \neq 0 .
\end{equation}
It follows from (\ref{e kNhN}) that
\begin{equation} \label{e hg-k2=f}
h_N (\lambda) g_N (\lambda) - k^2_N (\lambda) = f_N (\lambda)
\quad \mathrm{and} \quad h (\lambda) g (\lambda) - k^2
(\lambda) = f (\lambda).
\end{equation}
As above, the latter yields
       \begin{equation}\label{e M infZ}
M_\pm (\lambda ) = \frac {g(\pm \lambda)}{ k(\pm \lambda) \mp i
\sqrt{f(\pm \lambda)}} =  \frac {k(\pm \lambda) \pm i \sqrt{f(\pm
\lambda)}  }{ h (\pm \lambda)} .
      \end{equation}

$(a)$ Let $\ k(0) = 0$. Then (\ref{e f=0 gneq0}) and the first
equality in (\ref{e M infZ}) yield that (\ref{e asM Case1}) holds
with $C_1 = \prod_{j=1}^{\infty}  \xi_j (\prod_{j=1}^{\infty}
\mul_j \mur_j )^{-1/2} >0 $ (as above, the product converges due
to (\ref{e SumGap})).

$(b)$ Let $\ k(0) \neq 0$. Then (\ref{e hg-k2=f}) and (\ref{e f=0
gneq0}) yield $h(0) \neq 0$. Using the second representation of
$M_\pm (\lambda)$ from (\ref{e M infZ}), we get (\ref{e asM
Case2fz}) with the constants
 $C_2 = k(0)/h(0) \in \Rset$, $C_3 = |C_1/h(0)| >0$.

Theorem \ref{th_III.02} and Proposition \ref{p infReg}
complete the proof. 
\QQEEDD
\end{proof}

If the potential $q$ is periodic or finite(infinite)-zone and
$\inf \sigma(L) =0$, it is easy to show that $0$ \emph{is} a
critical point of $A$. So we have proved that $0$ is \emph{a
regular critical point} in these cases.

\section{Operators with nontrivial weights}\label{sec_IV}

In this section, we consider the $J$-self-adjoint operator $A$ of the type (\ref{I_02}) assuming
that $q\equiv 0$. In this case assumption (\ref{hyp_limpoint}) is fulfilled if and only if
$x\notin L^2(\Rset_\pm, |r(x)|dx)$. In the following
$\omega(\cdot)$ stands for $|r(\cdot)|$. Let us denote the
corresponding operator by
\begin{equation}\label{A_omega}
A_\omega:=-\frac{(\sgn x)}{\omega(x)}\frac{d^2}{dx^2},\qquad\qquad
\dom(L_\omega)=\D.
\end{equation}
Note that the operator $A_\omega$ is $J$-nonnegative. Hence,
by Proposition \ref{p RealS}, the spectrum of $A_\omega$ is real,
$\sigma(A_\omega)\subset \Rset$.


The main aim of this section is to prove 
Theorem \ref{th_IV.02}. But first we need two preparatory lemmas.

Consider the spectral problem
\begin{equation}\label{V_EvZ}
-y''(x)=\lambda x^\alpha y(x),\qquad x\geq0; \qquad y'(0)=0,
\end{equation}
with $\alpha>-1$. Denote by $z^{1/(2+\alpha)},\ z\in\Cset \setminus \Rset_+$, the
branch of the multifunction with cut along $\Rset_+$ such that $(-1)^{1/(2+\alpha)}=e^{i\pi/(2+\alpha)}$.
\begin{lem}[\cite{EvZ78}]\label{lem_IV_01}
Let $\alpha >-1$. Then the function
      \begin{equation}\label{IV_05}
            m_\alpha(\lambda):=C_{\nu}e^{i\pi\nu}\lambda^{-\nu}, \qquad \lambda\in\Cset_+; \quad\nu=\frac{1}{\alpha+2},
       \quad C_{\nu}:= \frac{\Gamma(1+\nu)}{\nu^{2\nu}\Gamma(1-\nu)}
            ,
      \end{equation}
is the Titchmarsh-Weyl $m$-coefficient of the problem
(\ref{V_EvZ}).
Here $\Gamma(\cdot)$ is the classical $\Gamma$-function.
\end{lem}

This result was obtained 
in \cite{EvZ78} using an explicit form of the Weyl solution of
equation (\ref{V_EvZ})
(see \cite[Part III, equation 2.162 (1a)]{Kamke}). A different and simpler proof of Lemma \ref{lem_IV_01}
was given in \cite{Kos06} (but without computing $C_\nu$).

As a corollary of Lemma
\ref{lem_IV_01}, we obtain a simple proof of \cite[Theorem
2.7]{FN98}.
\begin{thm}[\cite{FN98}]
If $\omega(x)=|x|^\alpha$, \ $\alpha>-1$, then $A_\omega$ is similar to a self-adjoint operator in $L^2(\Rset, \omega(x)dx)$.
\end{thm}
\begin{proof}
The operator $A_{|x|^\alpha}$ is $J$-self-adjoint since
$x\notin L^2(\Rset_\pm, |x|^{\alpha}dx)$. By Lemma
\ref{lem_IV_01}, we have $ M_+(\cdot)= -M_-(-\cdot)=
m_{\alpha}(\cdot)$. Hence,
\begin{equation}\label{III_2_01}
\frac{M_+(\lambda)+M_-(\lambda)}{M_+(\lambda)-M_-(\lambda)}=
\frac{1+\exp{\{i\pi\nu\}}}{1-\exp{\{i\pi\nu\}}},\qquad
\lambda\in\Cset_+.
\end{equation}
By Theorem \ref{th_III.1}, $A_{|x|^\alpha}$ is similar to a
self-adjoint operator.\QQEEDD
\end{proof}

\begin{lem}\label{lem_IV_02}
Let $\alpha>-1$ and let $p$ be a
positive function satisfying (\ref{IV_03}) on $\Rset_+$ with
$\alpha_+=\alpha$ and certain $c_+>0$. Let $m_+(\cdot)$ be the
Titchmarsh-Weyl $m$-coefficient of the problem
\begin{equation}\label{IV_11}
-\frac{d^2 y(x)}{dx^2}=\lambda p(x)|x|^{\alpha}y(x),\qquad
x>0;\qquad y'(0)=0.
\end{equation}
Then
\[  
m_{+}(\lambda)=
C_{\nu}e^{i\pi\nu}(c_+\lambda)^{-\nu}(1 +o(1)), \qquad \lambda\to
0;\qquad \nu=\frac1{2+\alpha}.
 \]   
\end{lem}
\begin{proof}
Without loss of generality it can be assumed that $c_+=1$.

Let $f_p(x, \lambda)$ denote the Weyl
solution of (\ref{IV_11}) for $\lambda\in\Cset_+$ (we will write $f_1(x, \lambda)$ if
$p\equiv 1$).
It is known (see \cite[Part III, equation 2.162 (1a)]{Kamke}) that the general
solution of (\ref{V_EvZ}) is
\[  
y(x,\lambda)=c_1\sqrt{x}H_\nu^{(1)}(2\nu\sqrt{\lambda}x^{1/2\nu})+c_2\sqrt{x}H_\nu^{(2)}(2\nu\sqrt{\lambda}x^{1/2\nu}),
\]  
where $c_j\in\Cset$ and $H_\nu^{(j)}(\cdot)$ are the Hankel
functions (see \cite[\Sec 3.6]{Watson}). Moreover (see \cite[\Sec
7.2]{Watson}), if
$-\pi+\delta\leq \arg z\leq\pi-\delta$, $\delta>0$, then
\begin{eqnarray}
H_\nu^{(1)}(z)&=&\left(\frac{2}{\pi
z}\right)^{1/2}e^{i(z-\nu\pi/2-\pi/4)}(1+O(z^{-1})),\qquad |z|\to\infty;\label{IV_15}\\
H_\nu^{(2)}(z)&=&\left(\frac{2}{\pi
z}\right)^{1/2}e^{-i(z-\nu\pi/2-\pi/4)}(1+O(z^{-1})),\qquad
|z|\to\infty.\label{IV_16}
\end{eqnarray}
Note that (\ref{IV_15})-(\ref{IV_16}) implies  $f_1(x,\lambda)=\sqrt{x}H_\nu^{(1)}(2\nu\sqrt{\lambda}x^{1/2\nu})$ and
\begin{equation}\label{Hankel_det}
W(H_\nu^{(1)}(z), H_\nu^{(2)}(z))=-\frac{4i}{\pi z},
\end{equation}
where $W$ is the Wronskian $W(f,g)(x):=f(x)g'(x)-f'(x)g(x)$.

Let us consider the Green function
\[ 
G(x,t; \lambda)=\varphi_1(t, \lambda)f_1(x, \lambda)-\varphi_1(x,
\lambda)f_1(t, \lambda).
\] 
Here $\varphi_1=c_2(\lambda)\sqrt{x}H_\nu^{(2)}(2\nu\sqrt{\lambda}x^{1/2\nu})$ and $c_2$ is chosen such that $W(f_1,\varphi_1)\equiv 1$ (cf. (\ref{Hankel_det})).
Using (\ref{IV_16})--(\ref{Hankel_det}), after straightforward
calculations we obtain
\[
\left|G(x,t; \lambda)\frac{f_1(t; \lambda)}{f_1(x;
\lambda)}\right| \leq C_1 \left| \frac{\lambda^{-1/2}}{(
t^{\alpha/2}+1)} \left(
1-\exp\{i\lambda^{1/2}2\nu(|t|^{1/2\nu}-|x|^{1/2\nu})\}\right)\right|.
\]
Hence, for $0<x<t$,
\begin{equation}\label{IV_17}
\left|G(x,t; \lambda)\frac{f_1(t; \lambda)}{f_1(x;
\lambda)}\right| \leq 2C_2 \left|
\frac{\lambda^{-1/2}}{t^{\alpha/2}+1}\right|.
\end{equation}
Consider the following integral equation
\begin{equation}\label{IV_12}
y_p(x; \lambda)=f_1(x;
\lambda)+\lambda\int_x^{+\infty}|t|^{\alpha}(p(t)-1)G(x,t;
\lambda)y_p(t; \lambda)dt.
\end{equation}
Using a standard technique (see, for example, \cite[\Sec 3.1]{Mar77}), one can show that (\ref{IV_17}) and (\ref{IV_03}) imply that the solution of (\ref{IV_12}) exists and is the Weyl solution of (\ref{IV_11}).
Denoting $y_p(x, \lambda)=f_1(x,
\lambda)\widetilde{y_p}(x,\lambda)$ in (\ref{IV_12}), one gets
\begin{equation}\label{IV_18}
\widetilde{y_p}(x;
\lambda)=1+\lambda\int_x^{+\infty}|t|^{\alpha}(p(t)-1)G(x,t;
\lambda)\frac{f_1(t; \lambda)}{f_1(x; \lambda)}\widetilde{y_p}(t,
\lambda)dt.
\end{equation}
Combining (\ref{IV_18}) with (\ref{IV_17}) and (\ref{IV_03}), we arrive
at
\begin{equation}\label{IV_19}
f_p(0, \lambda)=f_1(0, \lambda)(1+o(1)),\qquad \lambda\to 0.
\end{equation}
Analogously one obtains
\begin{equation}\label{IV_20}
f_p'(0, \lambda)=f_1'(0, \lambda)(1+o(1)),\qquad \lambda\to 0.
\end{equation}
Combining (\ref{IV_19}), (\ref{IV_20}) with (\ref{def_wf}), we obtain 
\[
m_+(\lambda)=-\frac{f_p(0, \lambda)}{f_p'(0,
\lambda)}=-\frac{f_1(0, \lambda)}{f_1'(0,
\lambda)}(1+o(1))=m_\alpha(\lambda)(1+o(1)),\qquad \lambda\to 0.
\]
\QQEEDD
\end{proof}

\begin{proof}[Proof of Theorem \ref{th_IV.02}]
$(i)$ By (\ref{IV_03}) and Lemma \ref{lem_IV_02}, we obtain
\begin{eqnarray*}
M_+(\lambda)&=&m_+(\lambda)=C_{\nu_+}e^{i\pi\nu_+}(c_+\lambda)^{\nu_+}(1
+o(1)),\\
M_-(\lambda)&=&-m_-(-\lambda)=C_{\nu_-}(c_-\lambda)^{\nu_-}(1 +
o(1)),\qquad\lambda\to 0,\quad\lambda\in\Cset_+,
\end{eqnarray*}
where $\nu_\pm=1/(2+\alpha_\pm)$ and $c_\pm>0$. Therefore,
\[
\left|\frac{M+(\lambda)+M_-(\lambda)}{M+(\lambda)-M_-(\lambda)}\right|=
\left|\frac{C_{\nu_+}e^{i\pi\nu_+}(c_+\lambda)^{\nu_+}
-C_{\nu_-}(c_-\lambda)^{\nu_-}}{C_{\nu_+}e^{i\pi\nu_+}(c_+\lambda)^{\nu_+}+C_{\nu_-}(c_-\lambda)^{\nu_-}}
(1+o(1))\right|, \quad \lambda \to 0.
\]
Hence
$ \left( M+(\lambda)+M_-(\lambda)\right) \left(M+(\lambda)-M_-(\lambda)\right)^{-1}$ is bounded in a neighborhood of $0$. Thus, by Theorem \ref{th_III.02}, $0$ is
not a singular critical point of $A_\omega$.

$(ii)$ Condition (\ref{IV_03}) implies that $\omega\notin
L^1(\Rset_\pm)$. Hence $\ker A_\omega=\{0\}$. Combining (i) with
Propositions \ref{p infReg} and \ref{p cr=sim}, we obtain
the similarity of $A_\omega$ to a self-adjoint
operator. \QQEEDD
\end{proof}

 \begin{rem}\label{rem_th_FSh}
Using another approach, M.M.~Faddeev and R.G.~Shterenberg proved the similarity of $A_\omega$ to a self-adjoint operator  under additional rather strong  assumptions on the weight $\omega$ (see \cite[Theorem 7]{FSh2}).
We avoid these difficulties using the spectral theory of
$J$-nonnegative operators.
\end{rem}

\ack{ The authors thank L. Oridoroga
for useful discussions.
AK gratefully acknowledges support from the Junior Research
Fellowship Program of the Erwin Schr\"odinger Institute for
Mathematical Physics. IK would like to thank to the University of
Calgary for support from the Postdoctoral Fellowship Program.}

\quad
\\
Illya M. Karabash, \\
Department of Mathematics and Statistics,
 University of Calgary, \\2500 University Drive NW,
 Calgary T2N 1N4, Alberta, CANADA \\
 and \\
Institute of Applied Mathematics and Mechanics, NAS of Ukraine,\\
R. Luxemburg str., 74, Donetsk 83114, UKRAINE\\
\emph{e-mail:} karabashi$@$yahoo.com, karabashi$@$mail.ru\\
\\
Aleksey S. Kostenko, \\
Institute of Applied Mathematics and Mechanics, NAS of Ukraine,\\
R. Luxemburg str., 74, Donetsk 83114, UKRAINE\\
\emph{e-mail:} duzer80$@$mail.ru; \ duzer80$@$gmail.com\\
\\
Mark M. Malamud, \\
Institute of Applied Mathematics and Mechanics, NAS of Ukraine,\\
R. Luxemburg str., 74, Donetsk 83114, UKRAINE\\
\emph{e-mail:}\ mmm$@$telenet.dn.ua

\end{document}